\documentclass[a4paper,11pt,oneside]{article}
\usepackage[english]{babel}
\usepackage{graphicx}
\usepackage[latin1]{inputenc}
\usepackage{amssymb,amsmath,textcomp,amsthm, mathtools}
\usepackage{amsfonts}
\usepackage{setspace}
\usepackage{multicol, multirow}
\usepackage{color}

\usepackage{algorithm}
\usepackage{algo_olm}

\usepackage{pgfplots}

\newtheorem{thm}{Theorem}[section]

\newtheorem{definition}[thm]{Definition}

\newtheorem{property}[thm]{Property}

\newtheorem{remark}[thm]{Remark}

\newlength\figureheight
\newlength\figurewidth

\begin{document}
\title{A nonintrusive method to approximate linear systems with nonlinear parameter dependence}
\maketitle

\begin{center}
{Fabien Casenave$^{1}$, Alexandre Ern$^{1}$, Tony Leli\`{e}vre$^{1,2}$, and Guillaume Sylvand$^{3}$}\\
\end{center}
\begin{center}
\end{center}
\begin{center}
$^1$ Universit\'{e} Paris-Est, CERMICS (ENPC), 6-8 Avenue Blaise Pascal, Cit\'{e} Descartes ,  F-77455 Marne-la-Vallé\'{e}, France
\end{center}
\begin{center}
$^2$ INRIA Rocquencourt, MICMAC Team-Project, Domaine de Voluceau, B.P. 105, 78153 Le Chesnay Cedex, France 
\end{center}
\begin{center}
$^3$ EADS-IW, 18 rue Marius Terce, 31300 Toulouse, France
\end{center}

\begin{abstract}
We consider a family of linear systems $A_\mu \alpha=C$ with system matrix $A_\mu$ depending on a parameter $\mu$ and for simplicity
parameter-independent right-hand side $C$. These linear systems typically result from the finite-dimensional approximation of
a parameter-dependent boundary-value problem. We derive a procedure based on the Empirical Interpolation Method to obtain a separated
representation of the system matrix in the form $A_\mu\approx\sum_{m}\beta_m(\mu)A_{\mu_m}$ for some selected values of the
parameter. 
Such a separated representation is in particular useful in the Reduced Basis Method. The procedure is
called nonintrusive since it only requires to access the matrices $A_{\mu_m}$. As such, it offers a crucial advantage over
existing approaches that instead derive separated representations requiring to enter the code at the level of assembly.
Numerical examples illustrate the performance of our new procedure on a simple one-dimensional boundary-value problem and
on three-dimensional acoustic scattering problems solved by a boundary element method.
\end{abstract}

\section{Introduction}
In industrial projects, decisions are often taken after a series of complex computations using computer codes of
various origins. To simplify the overall computation, surrogate models can be used to replace some parts of the
computation. Some of these surrogates are constructed using only a series of input/output couples.
With some hypotheses on the input, confidence intervals can be derived, see e.g. \cite{krig} for the kriging method.
When additional knowledge on the underlying mathematical formulation is available, model reduction methods can be used.
For instance, the Reduced Basis Method (RBM) enables fast resolutions on a basis of precomputed solutions, rather than on a finite
element basis (see \cite{RB} for a detailed presentation and \cite{RBconv} for some convergence results).
We consider a family of linear systems $A_\mu\alpha=C$ of order $n$, where $n$ is large.
For simplicity, we assume that the right-hand side is independent of the parameter $\mu$.

The RBM consists first in an offline stage, where a reduced basis of $\hat{n}$ functions $u_j$, $1\leq j\leq \hat{n}$ are
computed using a greedy algorithm. These functions are solution of the original problem for some values
$\left(\mu_j\right)_{1\leq j\leq\hat{n}}$ of the parameter $\mu$, which are selected using a greedy algorithm. The functions
$u_j$ thus write $u_j(x)=\sum_{i=1}^n \alpha_i(\mu_j)\theta_i(x)$, where $\left(\theta_i\right)_{1\leq i\leq n}$ is the finite
element basis and the vector $\alpha(\mu_j)=\left(\alpha_i(\mu_j)\right)_{1\leq i\leq n}$ is such that $A_{\mu_j}\alpha(\mu_j)=C$.
In practice, the dimension of the reduced basis is much smaller than the dimension of the finite element basis: $\hat{n}\ll n$.
Denote by $U$ the rectangular matrix of size $n\times\hat{n}$ such that $(U)_{i,j}=\alpha_i(\mu_j)$.
Second, in the online stage, for a given value of the parameter $\mu$, a reduced problem is constructed as
$\hat{A}_\mu \hat{\alpha}(\mu)=\hat{C}$, where $\hat{A}_\mu=U^t A_\mu U$ and $\hat{C}=U^t C$. Solving this reduced problem for a
certain value of $\mu$ leads to the approximate solution $\hat{u}_\mu(x)=\sum_{j=1}^{\hat{n}}\hat{\alpha}_j(\mu)u_j(x)$. 

To efficiently construct the online problems, a separated representation (also known as an {\em affine decomposition} in the RBM
literature) of the matrix assembled by the code is needed in the form
\begin{equation}
\label{eq:cmame:affine_dep}
A_\mu \approx \sum_{m=1}^d\gamma_m(\mu)A_m, 
\end{equation}
so that
\begin{equation}
\label{eq:non-intrusive_RB}
\hat{A}_\mu\approx \sum_{m=1}^{d} \gamma_m(\mu) U^t A_m U,
\end{equation}
where the matrices $U^t A_m U$ are of small size $\hat{n}\times\hat{n}$ and can be precomputed during the offline stage.
The separated representation \eqref{eq:cmame:affine_dep} thus enables online problems to be constructed in complexity independent of $n$, as long as the functions
$\mu\mapsto \gamma_m(\mu)$ are also computed in complexity independent of $n$.
Standard techniques (see \cite{maday}) to obtain the separated representation \eqref{eq:cmame:affine_dep}
require in general nontrivial modifications of the assembling routines of the computational code in order to access
separately various terms of the variational formulation at hand (See Remark~\ref{remark_cmame} below for more details).

The present work provides a step forward in this context, since a procedure that yields a separated representation of $A_\mu$
in the form
\begin{equation}
\label{eq:cmame:goal}
A_\mu\approx\sum_{m=1}^d \beta_m(\mu)A_{\mu_m}
\end{equation}
is derived, where $\left(\mu_m\right)_{1\leq m\leq d}$ are some selected values of the parameter. Since the separated representation
\eqref{eq:cmame:goal} only uses the complete system matrix at the selected parameter values, this representation
requires no implementation effort in the assembly routines of the computational code under the (mild) assumptions that
we can indeed access the system matrix $A_\mu$ and that we can identify the functional dependencies on $\mu$ in the variational
formulation under consideration (see below for more details). For this reason, the procedure is called nonintrusive.

In Section~\ref{sec:approxpb}, we present the approximation problems investigated in this work, first a simple introductory
example and then problems with a more complex parameter dependence.
In Section~\ref{sec:EIM}, we briefly recall the Empirical Interpolation Method.
In Section~\ref{sec:non_intrusive}, we present our nonintrusive procedure for the introductory example and
test it on a one-dimensional boundary-value problem. The procedure is extended to more complex parameter dependence in Section~\ref{sec:applicacous}
where it is also applied to two three-dimensional scattering problems. Some conclusions are drawn in Section~\ref{sec:otherapprox}
where, in particular, we observe that
our procedure can be extended to the approximation of other quantities.

\section{The approximation problem}
\label{sec:approxpb}
We first present an introductory example.
Let $\mathcal{V}$ be a Hilbert space and
consider the following weak formulation: Find $u\in\mathcal{V}$ such that for all $u^t\in\mathcal{V}$,
\begin{equation}
\label{eq:varf}
\int_{\Omega} g(\mu,x){\nabla}u(x)\cdot {\nabla}u^t(x)dx+\int_{\Omega} \mu u(x)u^t(x)dx=b(u^t),
\end{equation}
where $\Omega$ is the domain of computation, $\mu$ a parameter belonging to a given parameter set $\mathcal{P}$, $g(\mu,x)$
a given function defined on $\mathcal{P}\times\Omega$ and $b$ a bounded linear form on $\mathcal{V}$.
Consider now a conforming $n$-dimensional approximation of the space $\mathcal{V}$ denoted by $\mathcal{V}_h$
(the subscript $h$ refers to an underlying mesh), and a basis of $\mathcal{V}_h$ denoted by $(\theta_i)_{1\leq i\leq n}$.
The finite element approximation of \eqref{eq:varf} requires the computation of the matrix $A_\mu$ of size $n\times n$
with entries
\begin{equation}
\label{eq:matrix}
\left(A_{\mu}\right)_{i,j}:=
\begin{pmatrix}
\displaystyle\int_{\Omega} g(\mu,x){\nabla}\theta_j(x)\cdot{\nabla}\theta_i(x)dx + \mu\int_{\Omega}\theta_j(x)\theta_i(x)dx
\end{pmatrix}_{i,j}.
\end{equation}
The notation $A_\mu$ is adopted to stress the fact that the matrix $A_\mu$ depends on the value of the parameter $\mu$.
The problem solved by the computational code is
\begin{equation}
\label{eq:matrix_pb}
A_\mu \alpha = C,
\end{equation}
where $(C)_i=b(\theta_i)$ for all $1\leq i\leq n$, and where an approximation of the solution $u$ to \eqref{eq:varf} is obtained in
the form $u(x)\approx \sum_{i=1}^n \alpha_i \theta_i(x)$.

Let
\begin{equation}
\label{eq:matrix1}
\left(A^1_{\mu}\right)_{i,j}:=
\begin{pmatrix}
\displaystyle\int_{\Omega} g(\mu,x){\nabla}\theta_j(x)\cdot{\nabla}\theta_i(x)dx
\end{pmatrix}_{i,j}\textnormal{~and~}
\left(A^0\right)_{i,j}:=
\begin{pmatrix}
\displaystyle\int_{\Omega} \theta_j(x)\theta_i(x)dx
\end{pmatrix}_{i,j}
\end{equation}
so that
\begin{equation}
\label{eq:cmame:simplexmat}
A_\mu=A^1_\mu+\mu A^0. 
\end{equation}

\begin{definition}[Intrusivity]
\label{defintrusivity}
A procedure leading to a separated representation of $A_\mu$ in the general form~\eqref{eq:cmame:affine_dep} is called
\begin{itemize}
 \item intrusive if it requires to implement new integral terms,
 \item weakly-intrusive if it only requires to precompute independently $A^1_{\mu}$ for some values of $\mu$ and $A_0$,
 \item nonintrusive if it only requires to precompute $A_\mu$ for some values of $\mu$.
\end{itemize}
\end{definition}

The term ``weakly-intrusive'' comes from the fact that the user has to enter the routines of the code and to insert switches at the
right places to save the terms in $A^1_{\mu}$ independently from the terms in $A^0$.
In the context of industrial codes, this is not always possible.
The notion of nonintrusivity in Definition \ref{defintrusivity} is different from the notion of black-box,
which requires only the computation of input / output couples.
Our purpose is to develop a nonintrusive procedure leading to the separated representation \eqref{eq:cmame:goal} of $A_{\mu}$.

The above example can be generalized to a class of engineering problems requiring to compute a large, parameter-dependent
matrix $A_\mu$ for many values of the parameter $\mu$ where $A_\mu$ is of the form
\begin{equation}
\label{eq:generalsetting}
A_\mu=\sum_{\varrho=1}^R A_\mu^\varrho+\sum_{s=1}^S \psi_s(\mu) A^{s},
\end{equation}
where $A_\mu^\varrho$ are matrices that require to integrate some functions $g^\varrho(\mu,x)$ over $\Omega$, $\psi^s$ are given functions
of $\mu$ and $A^{s}$ are $\mu$-independent matrices resulting from some integration over $\Omega$. The introductory example
corresponds to $R=1$, $S=1$, and $\psi_1(\mu)=\mu$.
To simplify the presentation of the main ideas,
we consider the setting of \eqref{eq:cmame:simplexmat} in Sections \ref{sec:EIM} and \ref{sec:non_intrusive}
and return to the more general setting of \eqref{eq:generalsetting} in Section \ref{sec:applicacous}.

\section{Empirical Interpolation Method}
\label{sec:EIM}
The Empirical Interpolation Method (EIM) is a procedure to approximate two-variable functions.
In particular, it can be used to approximate the two-variable function $g(\mu,x)$, for all $\mu\in\mathcal{P}$
and all $x\in\Omega$. Denote by ${\rm EIM}^g$ this particular procedure.
${\rm EIM}^g$ leads to an interpolation operator $I^g_{d^g}$ such that
\begin{equation}
\label{eq:EIMapprox}
\left(I^g_{d^g} g\right)(\mu,x)\approx g(\mu,x), \qquad\forall\mu\in\mathcal{P},~\forall x\in\Omega,
\end{equation}
where $d^g$ is the number of interpolation points (called {\em magic points} in the context of BRM, see \cite{maday}).
${\rm EIM}^g$ is composed of two stages: (i) an offline stage, where a matrix $B^g$ of size $d^g\times d^g$, a set of $d^g$
$x$-dependent basis functions $\{q^g_k\}_{1\leq k\leq d^g}$,
a set of $d^g$ points $\{x_k\}_{1\leq k\leq d^g}$ in $\Omega$, and a set a $d^g$ parameter values
$\{\mu_k\}_{1\leq k\leq d^g}$ in $\mathcal{P}$ are constructed, (ii) an online stage,
where the quantities computed in the offline stage are used to carry out the approximation~\eqref{eq:EIMapprox} (see
Section~\ref{sec:cmame:practimpl} for more details on the offline / online stages for the whole procedure).

The offline stage of ${\rm EIM}^g$ is detailed in Algorithm~\ref{algo1}. In the loop on $k$ in Algorithm~\ref{algo1}, the residual
operator $\delta^g_k$ is defined by $\delta^g_k={\rm Id}-{I}^g_{k}$, where the interpolation operator ${I}^g_{k}$ is such that
\begin{equation}
\label{eq:onlinea1}
\left(I^g_{k} g\right)(\mu,x) := \sum_{m=1}^{k} \lambda^g_m(\mu) q^g_{m}(x),
\end{equation}
and for a given $\mu\in\mathcal{P}$, the $\lambda_m^g(\mu)$'s are defined by
\begin{equation}
\label{eq:onlinea1pb}
\sum_{m=1}^{k}B^g_{l,m}{\lambda}^g_m(\mu)=g(\mu,x^g_l), \qquad \forall 1\leq l\leq k.
\end{equation}
After $d^g$ iterations, the interpolation formula~\eqref{eq:onlinea1} leads to the following approximation for $A_\mu$:
\begin{equation}
\label{eq:intrusive1}
A_\mu\approx \sum_{m=1}^{d^g}\lambda^g_m(\mu)M_m+\mu A^0,
\end{equation}
where $\left(M_m\right)_{i,j}=\int_{\Omega} q^g_{m}(x)\vec{\nabla}\theta_j(x)\cdot\vec{\nabla}\theta_i(x)dx$. This representation
of $A_\mu$ is of the form~\eqref{eq:cmame:affine_dep}.

\begin{algorithm}[h!]
	\caption{Offline stage of ${\rm EIM}^g$}
	\label{algo1}
	\begin{algorithmic}[1]
        \STATE {Choose $d^g>1$}
        \hfill \COMMENT{Number of interpolation points}
	\STATE {Set $k:=1$}
        \STATE {Compute $\displaystyle \mu^g_1:=\underset{\mu\in \mathcal{P}}{\textnormal{argmax}}\|g(\mu,\cdot)\|_{L^\infty\left(\Omega\right)}$}
        \STATE {Compute $\displaystyle x^g_1:=\underset{x\in\Omega}{\textnormal{argmax}}|g(\mu^g_1,x)|$}
        \hfill \COMMENT{First interpolation point}
	\STATE {Set $\displaystyle q^g_1(\cdot):=\frac{g(\mu^g_1,\cdot)}{g(\mu^g_1,x^g_1)}$}
       \hfill \COMMENT{First basis function}
        \STATE {Set $B^g_{1,1}:=1$}
       \hfill \COMMENT{Initialize $B^g$ matrix}
        \WHILE {$k\leq d^g$} 
		\STATE Compute $\displaystyle \mu^g_{k+1}:=\underset{\mu\in \mathcal{P}}{\textnormal{argmax}}\|(\delta^g_k g)(\mu,\cdot)\|_{L^\infty\left(\Omega\right)}$
                \STATE Compute $\displaystyle x^g_{k+1}:=\underset{x\in\Omega}{\textnormal{argmax}}|(\delta^g_k g)(\mu^g_{k+1},x)|$
                \hfill \COMMENT{$(k+1)$-th interpolation point}  
                \STATE Set $\displaystyle q^g_{k+1}(\cdot):=\frac{(\delta^g_k g)(\mu^g_{k+1},\cdot)}{(\delta^g_k g)(\mu^g_{k+1},x^g_{k+1})}$
                \hfill \COMMENT{$(k+1)$-th basis function}
                \STATE Set $\displaystyle B^g_{i,{k+1}}:=q^g_{k+1}(x^g_i)$, for all $1\leq i\leq {k+1}$
                \hfill \COMMENT{Increment matrix $B^g$}
                \STATE $k\leftarrow k+1$
                \hfill \COMMENT{Increment the size of the interpolation}
	\ENDWHILE
\end{algorithmic}
\end{algorithm}

\begin{property}[Interpolation]
\label{propint1}
$\forall x\in\Omega,~\forall~1\leq m\leq d^g,~\left(I^g_{d^g} g\right)(\mu^g_m,x) = g(\mu^g_m,x)$.
\end{property}
\begin{proof}
See \cite[Lemma 1]{maday}.
\end{proof}
Property~\ref{propint1} means that, at the parameter values $\left(\mu^g_k\right)_{1\leq k\leq d^g}$ selected by ${\rm EIM}^g$, the
approximation~\eqref{eq:intrusive1} is exact since $\sum_{m=1}^{d^g}\lambda_m^g(\mu_k^g)M_m=A^1_{\mu_k^g}$ for all $1\leq k\leq d^g$.
Since ${\rm Vect}_{1\leq k\leq d^g}\left(q^g_k(x)\right)={\rm Vect}_{1\leq k\leq d^g}\left(g(\mu^g_k,x)\right)$ holds in
Algorithm \ref{algo1}, the functions $q^g_k(x)$ can be expressed in terms of the functions $g(\mu^g_k,x)$ in the following
form: there exist $\gamma_{l,k}$, $1\leq l\leq k \leq d^g$ such that
$q^g_k(x)=\sum_{l=1}^{d^g} \gamma_{l,k}g(\mu^g_l,x)$, for all $1\leq k\leq d^g$.
Letting $(\lambda^g_m(\mu))_{1\leq m\leq d^g}$ solve \eqref{eq:onlinea1pb} for $k=d^g$, we obtain after exchanging the summations
\begin{equation}
(I^g_{d^g} g)(\mu,x)=\sum_{m=1}^{d^g}\left(\sum_{l=1}^{d^g}\gamma_{m,l}\lambda^g_l(\mu)\right)g(\mu^g_m,x).
\end{equation}
Define $\displaystyle \eta^g_m(\mu):=\sum_{l=1}^{d^g}\gamma_{m,l}\lambda^g_l(\mu)$.
The following property then holds:
\begin{property}[Weak-intrusivity]
${\rm EIM}^g$ leads to a weakly-intrusive procedure, since the resulting approximation of $A_\mu$ can be written
\begin{equation}
\label{eq:weakly_int}
A_\mu \approx \sum_{m=1}^{d^g}\eta^g_m(\mu)A^1_{\mu^g_m}+\mu A^0.
\end{equation}
\end{property}

\begin{remark}[Comparison with the standard EIM procedure in the RBM literature]
\label{remark_cmame}
If the considered variational formulation contains only one term, the above procedure was already proposed in the RBM literature
as a nonintrusive method to obtain a separated representation of the linear system under consideration, see~\cite{maday}.
For instance in~\eqref{eq:weakly_int}, if $A^0=0$, then $A^1_{\mu^g_m}=A_{\mu^g_m}$.
In the general setting of~\eqref{eq:generalsetting}, this corresponds to $R=1$ and $S=0$.
In any other case, the classical EIM needs to access independently matrices associated to each term of the variational formulation
and thus cannot deliver a separated representation solely based on the $A_\mu$ matrices.
\end{remark}

\section{The nonintrusive procedure}
\label{sec:non_intrusive}
\subsection{Description of the procedure}
\label{sec:descriptionproc}

Denote by $G^g(\mu)$ the vector-valued function with $d^g$ components such that $G^g_m(\mu)=g(\mu,x^g_m)$, for all $1\leq m\leq d^g$.
Then, from~\eqref{eq:onlinea1pb}, $\lambda^g(\mu)=(\lambda^g_m(\mu))_{1\leq m\leq d^g}$ can be concisely written as $\lambda^g(\mu)=\left(B^g\right)^{-1}G^g(\mu)$.
Notice that the computation of $\lambda^g(\mu)$ only requires the matrix $B^g$ and the set of points $\{x_m^g\}_{1\leq m\leq d^g}$.
Let $\left(z_p(\mu)\right)_{1\leq p\leq d_{\rm max}}$ with $d_{\rm max}:=d^g+1$, be such that
\begin{equation}
\label{eq:defX}
z_p(\mu):=\left\{
\begin{alignedat}{2}
&\lambda^g_p(\mu)&\quad &1\leq p\leq d^g,\\
&\mu&\quad &p=d^g+1.
\end{alignedat}\right.
\end{equation}
Recalling the notation $(M_m)_{i,j}:=\int_{\Omega} q^g_{m}(x)\vec{\nabla}\theta_j(x)\cdot\vec{\nabla}\theta_i(x)dx$
for all $1\leq m\leq d^g$, we infer from \eqref{eq:intrusive1} that
\begin{equation}
\label{eq:linform}
A_\mu\approx \sum_{m=1}^{d^g}\lambda^g_m(\mu)M_m+\mu A^0=\sum_{p=1}^{d_{\rm max}}z_p(\mu)T_p ,
\end{equation}
where the matrices
\begin{equation}
T_p:=\left\{
\begin{alignedat}{2}
&M_p&\quad &1\leq p\leq d^g,\\
&A^0&\quad &p=d^g+1=d_{\rm max},
\end{alignedat}\right.
\end{equation}
are independent of $\mu$. Note that $d_{\rm max}$ is the number of matrices to precompute and store when using the approximation
\eqref{eq:intrusive1}.

The key idea is now to apply a second EIM to approximate $z_p(\mu)$, where $z$ is seen as a function depending on the two
variables $p$ and $\mu$.
The EIM procedure to approximate $z_p(\mu)$ is denoted by ${\rm{EIM}}^z$ and its offline stage is detailed in Algorithm~\ref{algo3}.
The number of interpolation points is denoted by ${d}^z\leq {d}_{\rm max}$.
In the loop on $k$ in Algorithm~\ref{algo3}, the residual operator ${\delta}^z_k$ is defined by ${\delta}^z_k={\rm Id}-{I}^z_{k}$,
where
\begin{equation}
\label{eq:EIM2}
\left({I}^z_{k}z\right)_p(\mu):= \sum_{m=1}^{k} \beta^z_{m}(\mu) z_{p}(\mu^z_m),
\end{equation}
and
\begin{equation}
\label{eq:EIM2online}
\sum_{m=1}^{k}{B}^z_{m,l}\beta^z_{m}(\mu)={q}^z_l(\mu),\qquad 1\leq l\leq k.
\end{equation}
Owing to the interpolation property, there holds $(I^z_{d^z}z)_{{p}^z_k}(\mu)=z_{{p}^z_k}(\mu)$
for all $1\leq k\leq d^z$ and all $\mu\in \mathcal{P}$.
If $d^z=d_{\rm max}$, all the indices $p$ are selected in Algorithm \ref{algo3} and
$(I^z_{d_{\rm max}}z)_{p}(\mu)=z_p(\mu)$ for all $1\leq p\leq d_{\rm max}$ and all $\mu\in \mathcal{P}$.
Observe that we can stop ${\rm EIM}^z$ before ${d}^z=d_{\rm max}$ interpolation matrices have been computed,
see Sections~\ref{sec:cmame:num1} and~\ref{sec:application2} for some illustrations.

\begin{algorithm}[h!]
	\caption{Offline stage of ${\rm EIM}^z$}
	\label{algo3}
	\begin{algorithmic}[1]
        \STATE {Choose ${d}^z>1$}
        \hfill \COMMENT{Number of interpolation points}
	\STATE {Set $k:=1$}
        \STATE {Compute $\displaystyle {p}^z_1:=\underset{1\leq p\leq d^g+1}{\textnormal{argmax}}\|(z)_p(\cdot)\|_{L^\infty(\mathcal{P})}$}
        \STATE {Compute $\displaystyle {\mu}^z_1:=\underset{\mu\in \mathcal{P}}{\textnormal{argmax}}|(z)_{{p}^z_1}(\mu)|$}
        \hfill \COMMENT{First interpolation point}
	\STATE {Set $\displaystyle {q}^z_1(\cdot):=\frac{(z)_{{p}^z_1}(\cdot)}{(z)_{{p}^z_1}({\mu}^z_1)}$}
       \hfill \COMMENT{First basis function}
        \STATE {Set ${B}^z_{1,1}:=1$}
       \hfill \COMMENT{Initialize ${B}^z$ matrix}
        \WHILE {$k\leq {d}^z$} 
		\STATE Compute $\displaystyle {p}^z_{k+1}:=\underset{1\leq p\leq d^g+1}{\textnormal{argmax}}\|({\delta}^z_k z)_p(\cdot)\|_{L^\infty(\mathcal{P})}$,
                \STATE Compute $\displaystyle \mu^z_{k+1}:=\underset{\mu\in \mathcal{P}}{\textnormal{argmax}}|({\delta}^z_k z)_{{p}^z_{k+1}}(\mu)|$
                \hfill \COMMENT{$(k+1)$-th interpolation point}  
                \STATE Set $\displaystyle {q}^z_{k+1}(\cdot):=\frac{(\delta^z_k z)_{{p}^z_{k+1}}(\cdot)}{(\delta^z_k z)_{{p}^z_{k+1}}(\mu^z_{k+1})}$
                \hfill \COMMENT{$(k+1)$-th basis function}
                \STATE $\displaystyle {B}^z_{i,{k+1}}:={q}^z_{k+1}(\mu^z_i)$, for all $1\leq i\leq {k+1}$
                \hfill \COMMENT{Increment matrix ${B}^z$}
                \STATE $k\leftarrow k+1$
                \hfill \COMMENT{Increment the size of the interpolation}
	\ENDWHILE
\end{algorithmic}
\end{algorithm}
Injecting the approximation \eqref{eq:EIM2} with $k={d}^z$ into the right-hand side of \eqref{eq:linform} with $z_p(\mu)$
replaced by $(I_{d^z}^z z)_p(\mu)$ yields
\begin{equation}
\label{eq:non-intrusive}
A_\mu\approx \sum_{p=1}^{d_{\rm max}}T_p \sum_{m=1}^{{d}^z} \beta^z_{m}(\mu) z_{p}(\mu^z_m)
=\sum_{m=1}^{{d}^z} \beta^z_{m}(\mu)\sum_{p=1}^{d_{\rm max}}T_p z_{p}(\mu^z_m)
\approx \sum_{m=1}^{{d}^z} \beta^z_{m}(\mu) A_{\mu^z_m},
\end{equation}
where $\beta_m^z(\mu)$ is obtained from \eqref{eq:EIM2online}.
The right-hand side of \eqref{eq:non-intrusive} is the desired separated representation of $A_\mu$ that can be built
in a nonintrusive way.

\subsection{Practical implementation}
\label{sec:cmame:practimpl}

To compute the $L^\infty$-norms and determine the argmax in Algorithms~\ref{algo1} and~\ref{algo3}, it is convenient to consider
finite subsets of $\mathcal{P}$ and $\Omega$, denoted respectively by $\mathcal{P}_{\rm trial}$ and $\Omega_{\rm trial}$.
This becomes necessary when, for instance, the function $g(\mu,x)$ is not known analytically, but only for some elements of
$\mathcal{P}$ and $\Omega$. It seems natural to take for $\Omega_{\rm trial}$ the set of Gauss points on which the quadrature
formulae to compute the integrals in \eqref{eq:varf} are defined. However, this supposes to know and manipulate the set of the Gauss
points associated with the mesh.
Since the functions $q^g$ defined in Algorithm~\ref{algo1} are only used to construct the matrix $B^g$
and are not directly integrated with respect to $x$ to carry out the interpolation~\eqref{eq:non-intrusive}, it is possible
to write the procedure with any set $\Omega_{\rm trial}$ sampling the geometry.
Such an approach is considered in the numerical example of Section~\ref{sec:application2}.
More generally, the sets $\mathcal{P}_{\rm trial}$ and $\Omega_{\rm trial}$ should be fine enough to capture all the
phenomena, but not too fine to limit the overall computational cost. The numerical examples of Section~\ref{sec:applicacous} indicate
that high accuracy can be obtained with simple choices for $\mathcal{P}_{\rm trial}$ and $\Omega_{\rm trial}$.

In addition to the two sets $\mathcal{P}_{\rm trial}$ and $\Omega_{\rm trial}$,
the number of interpolation points $d^g$ and $d^z$ for each EIM have to be chosen.
The choice we made is to stop the two EIM's when respectively
$(\delta^g_k g)(\mu^g_{k+1},x^g_{k+1})$ and $(\delta^z_k z)_{{p}^z_{k+1}}(\mu^z_{k+1})$ have reached a prescribed threshold,
typically set at the level of the machine precision.

Finally, we specify the offline and online stages of our procedure when used within the RBM.
${\rm EIM}^g$ and the offline stage of ${\rm EIM}^z$ are part of the offline stage of the RBM. During the online stage of the RBM,
the reduced matrix is constructed as
\begin{equation}
\hat{A}_\mu\approx \sum_{m=1}^{{d}^z} \beta^z_{m}(\mu) U^t A_{\mu^z_m} U,
\end{equation}
so that only the online stage of ${\rm EIM}^z$ (i.e., the resolution of~\eqref{eq:EIM2online}) is needed.

\subsection{Illustration}

As a first illustration, we consider the following boundary-value problem:
\begin{equation}
\label{eq:simple1dBVP}
-\frac{d}{dx}\left(\exp(\mu x)\frac{du}{dx}(x)\right)+\mu u(x) =1\qquad\textnormal{in }\Omega:=(-3,3),\\
\end{equation}
with the following Dirichlet boundary condition $u(-3)=u(3)=0$.
The weak form reads: Find $u\in H^1_0(\Omega)$ such that for all $u^t\in H^1_0(\Omega)$,
\begin{equation}
\label{eq:varf1d}
a_\mu(u,u^t) = \int_{\Omega} u^t(x) dx,
\end{equation}
with
\begin{equation}
a_\mu(u,u^t) :=\int_{\Omega} \exp(\mu x)\frac{du}{dx}(x)\frac{du^t}{dx}(x)dx+\int_{\Omega} \mu u(x)u^t(x)dx.
\end{equation}

\begin{figure}[h!]
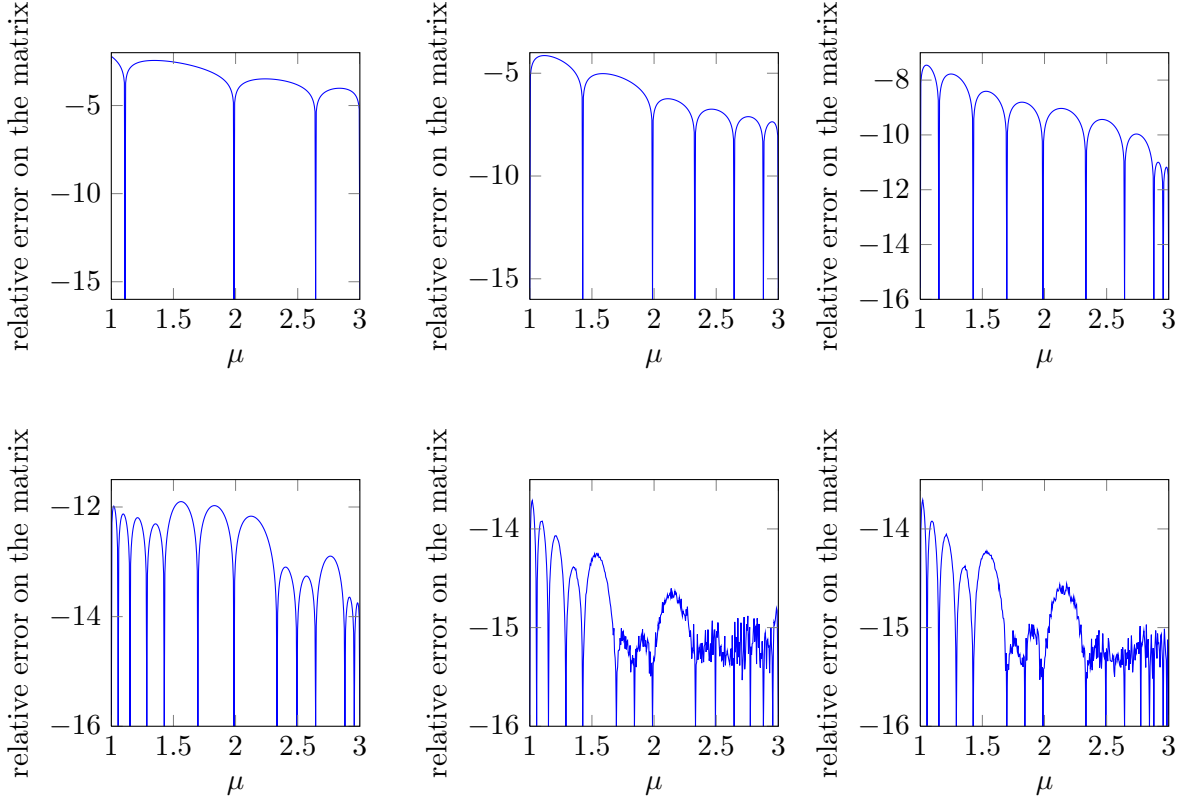

   \begin{minipage}[c]{.32\linewidth}
\include{erreurmatred3.tikz}
   \end{minipage} \hfill
   \begin{minipage}[c]{.32\linewidth}
\include{erreurmatred6.tikz}
   \end{minipage}
\begin{minipage}[c]{.32\linewidth}
\include{erreurmatred9.tikz}
   \end{minipage}
   \begin{minipage}[c]{.32\linewidth}
\include{erreurmatred12.tikz}
   \end{minipage} \hfill
   \begin{minipage}[c]{.32\linewidth}
\include{erreurmatred14.tikz}
   \end{minipage}
\begin{minipage}[c]{.32\linewidth}
\include{erreurmatred16.tikz}
   \end{minipage}
 \caption{Log${}_{10}$ of the relative Frobenius norm error on the matrix $A_\mu$ for $d^g=3, 6, 9, 12, 14$, $16$, and $d^z=d^g+1$.}
\label{fig:diffmat1d1}
\end{figure}

\begin{figure}[h!]
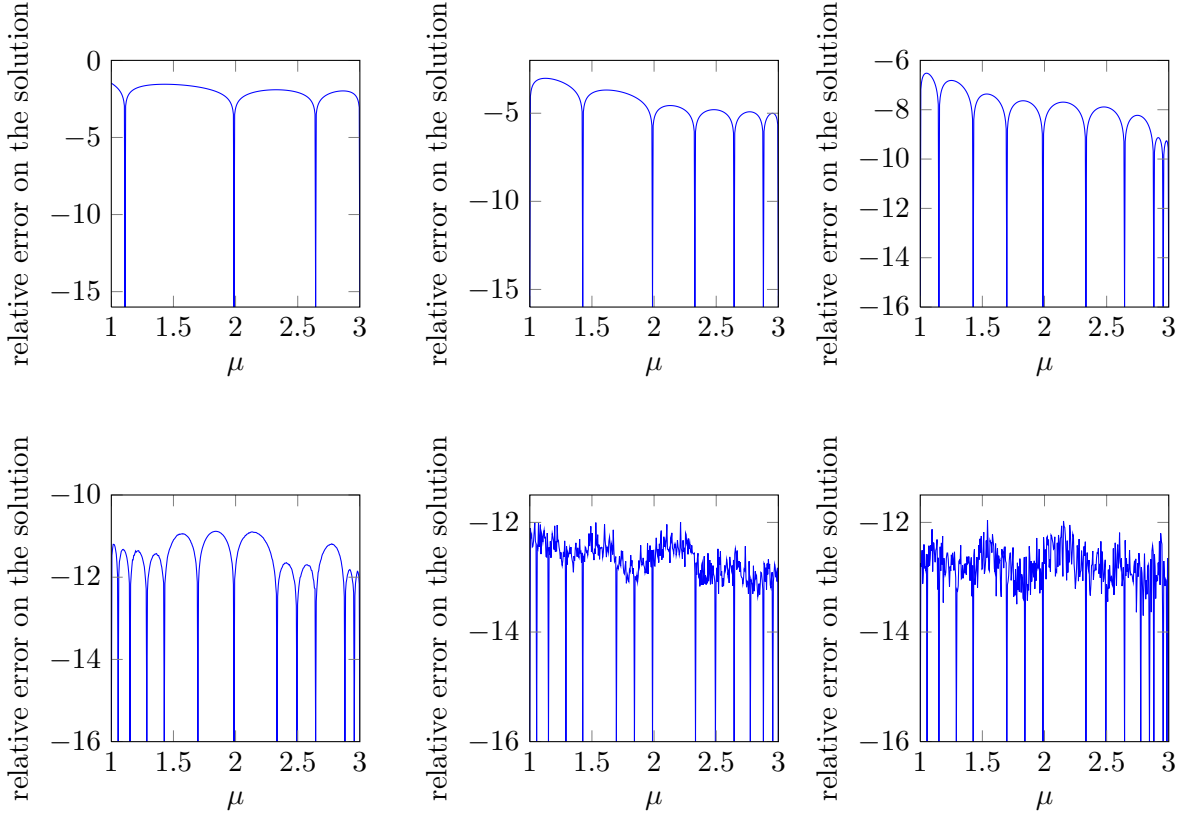

   \begin{minipage}[c]{.32\linewidth}
\include{erreursol3.tikz}
   \end{minipage} \hfill
   \begin{minipage}[c]{.32\linewidth}
\include{erreursol6.tikz}
   \end{minipage}
\begin{minipage}[c]{.32\linewidth}
\include{erreursol9.tikz}
   \end{minipage}
   \begin{minipage}[c]{.32\linewidth}
\include{erreursol12.tikz}
   \end{minipage} \hfill
   \begin{minipage}[c]{.32\linewidth}
\include{erreursol14.tikz}
   \end{minipage}
\begin{minipage}[c]{.32\linewidth}
\include{erreursol16.tikz}
   \end{minipage}
\caption{Log${}_{10}$ of the relative $L^2(\Omega)$-norm error on the solution for  $d^g=3, 6, 9, 12, 14$, $16$, and $d^z=d^g+1$.}
\label{fig:diffsol1d1}
\end{figure}

First-order continuous Lagrange finite elements are used, with a three-point quadrature formula in each mesh cell.
The mesh is uniform with $h_x=0.015$.
$\Omega_{\rm trial}$ is taken to be the set of Gauss points on the obtained mesh, and $\mathcal{P}_{\rm trial}=\{1, 1+h_\mu,
1+2h_\mu,... , 3\}$ with $h_\mu=0.005$.
To derive the separated approximation~\eqref{eq:non-intrusive}, ${\rm EIM}^g$ is first applied to $g(\mu,x):=\exp(\mu x)$.
Then, the vector-valued function $z(\mu)$ is constructed using \eqref{eq:defX}.
The quality of the whole procedure is measured, for various values of $d^g$ and $d^z=d^g+1$ using two error criteria:
(i) the relative Frobenius norm error on the matrix $A_\mu$ and (ii) the relative $L^2(\Omega)$-norm error on the solution,
see Figures~\ref{fig:diffmat1d1} and~\ref{fig:diffsol1d1}.

\setlength\figureheight{0.31\textwidth}
\setlength\figurewidth{0.31\textwidth}

We conclude from this first test case that the present method allows for a very good approximation of the matrix and the solution.

\section{Extension to more general parameter dependence}
\label{sec:applicacous}
The goal of this section is to show how to extend the nonintrusive procedure described in Section \ref{sec:non_intrusive} to more
complex parameter dependence. We illustrate the procedure on an industrial test case, namely a frequency-dependent three-dimensional
aeroacoustic scattering problem.

\subsection{Generalization of the nonintrusive procedure}

Recall the general form of the matrix $A_\mu$ to approximate:
\begin{equation}
A_\mu=\sum_{\varrho=1}^R A_\mu^\varrho+\sum_{s=1}^S \psi_s(\mu) A^{s},
\end{equation}
where $A_\mu^\varrho$ are matrices that require to integrate some functions $g^\varrho(\mu,x)$ over $\Omega$, $\psi^s$ are given functions
of $\mu$ and $A^{s}$ are $\mu$-independent matrices resulting from some integration over $\Omega$.
${\rm EIM}^g$ is applied independently to each $g^\varrho(\mu,x)$, for all $1\leq \varrho\leq R$, where the number of interpolation points,
respectively $(d^g)^\varrho$, may differ from one ${\rm EIM}^g$ to the other. These procedures lead to the construction of the functions
$(\lambda^g_m)^\varrho(\mu)$, for all $1\leq \varrho\leq R$, all $1\leq m\leq (d^g)^\varrho$, and all $\mu\in\mathcal{P}_{\rm trial}$,
using~\eqref{eq:onlinea1pb}. Then, define the functions $\left(z_p(\mu)\right)_{1\leq p\leq d_{\rm max}}$ with
$d_{\rm max}:=\sum_{\varrho=1}^{R}(d^g)^\varrho+S$ such that
\begin{equation}
z_p(\mu):=\left\{
\begin{alignedat}{2}
&(\lambda^g_m)^1(\mu),&\qquad &1\leq p\leq (d^g)^1,\quad  m=p,\\
&&\vdots&\\
&(\lambda^g_m)^R(\mu),&\qquad &1+\sum_{\varrho=1}^{R-1}(d^g)^\varrho\leq p\leq \sum_{\varrho=1}^{R}(d^g)^\varrho,\quad m=p-\sum_{\varrho=1}^{R-1}(d^g)^\varrho,\\
&\psi_1(\mu),&\quad &p=\sum_{\varrho=1}^{R}(d^g)^\varrho+1,\\
&&\vdots&\\
&\psi_S(\mu),&\quad &p=\sum_{\varrho=1}^{R}(d^g)^\varrho+S,
\end{alignedat}\right.
\end{equation}
and let ${\rm EIM}^z$ be applied to $z_p(\mu)$, with ${d}^z$ interpolation points, such that ${d}^z\leq {d}_{\rm max}=
\sum_{\varrho=1}^{R}(d^g)^\varrho+S$, to obtain an approximation of $A_\mu$ in the same form as~\eqref{eq:non-intrusive}.
Note that $d_{\rm max}$ is the number of matrices to precompute and store when using the approximation~\eqref{eq:intrusive1},
while the number of matrices to precompute and store when using~\eqref{eq:non-intrusive} is $d^z$; in our numerical examples
(see below), accurate representations of $A_\mu$ are already achieved for $d^z$ smaller than $d_{\rm max}$. Notice also
that in total, there are $\left(R+1\right)$ EIM procedures to be applied.

\subsection{Sound-hard scattering in the air at rest}
\label{sec:cmame:num1}

\begin{figure}[h!]
 \centering
\includegraphics[width=0.50\textwidth]{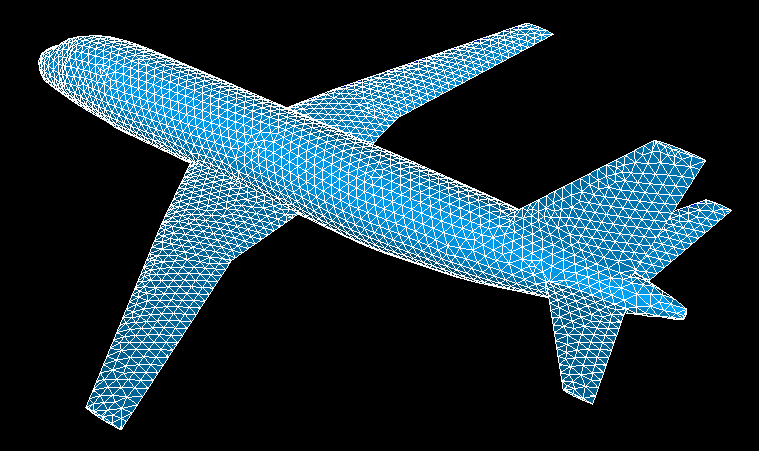}
\includegraphics[width=0.48\textwidth]{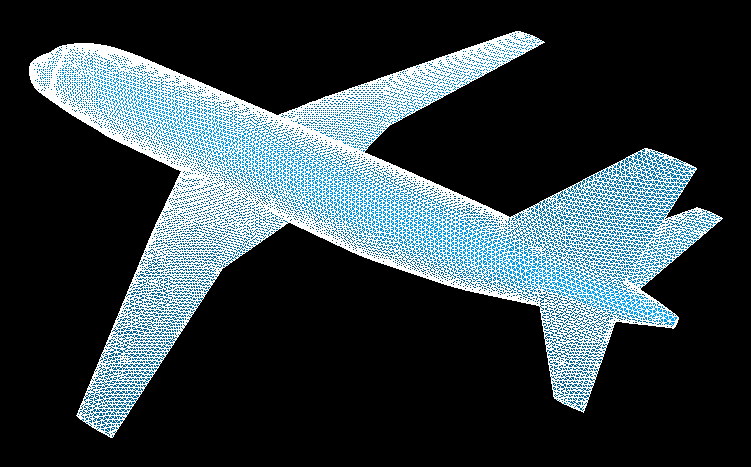}
 \caption{Airbus A319: Mesh~1 and Mesh~2.}
\label{mesh1}
\end{figure}

The problem of interest is the sound-hard scattering of an acoustic monopole source of wave number $\mu$ by an aircraft
(whose boundary is denoted by $\Gamma$)
in the air at rest, in the time-harmonic case. To simulate the noise created by one of the engines, the monopole is located under the
left wing of the plane. This is a classical Helmholtz exterior problem, for which one possible weak formulation is:
Find $u\in H^{\frac{1}{2}}\left(\Gamma\right)$ such that for all $u^t\in H^{\frac{1}{2}}\left(\Gamma\right)$,
\begin{equation}
\label{eq:hardsound}
a_\mu(u,u^t)=\int_{\Gamma}{f_{\rm inc}}_\mu(x) u^t(x) dx,
\end{equation}
where
\begin{equation}
\label{eq:biliform}
\begin{aligned}
a_\mu(u,u^t)&:=\frac{1}{4\pi}\int_{\Gamma}\int_{\Gamma}\frac{\exp\left(i\mu\left|x-y\right|\right)}{\left|x-y\right|}
\left(\overrightarrow{\rm curl}_{\Gamma}u(x)\cdot\overrightarrow{\rm curl}_{\Gamma}u^t(y)\right)dxdy\\
&\quad -\frac{\mu^2}{4\pi}\int_{\Gamma}\int_{\Gamma}\frac{\exp\left(i\mu\left|x-y\right|\right)}{\left|x-y\right|}
u(x)u^t(y)\left(\overrightarrow{n_x}\cdot\overrightarrow{n_y}\right)dxdy,
\end{aligned}
\end{equation}
where $\overrightarrow{\rm curl}_{\Gamma}$ denotes the surfacic curl on $\Gamma$, $\overrightarrow{n_x}$ the unit normal vector on
$\Gamma$ pointing towards the medium of propagation,
and ${f_{\rm inc}}_\mu$ is the incident acoustic field created by the source.
We refer to \cite[Section 3.4]{nedelecbook} for details on the derivation of \eqref{eq:hardsound}, and justifications on the
well-posedness of the integral in~\eqref{eq:biliform}.
The parameter of interest is the wave number $\mu$ of the acoustic monopole source.
The Boundary Element Method (BEM) is used to approximate problem~\eqref{eq:hardsound}. This leads to a dense
$\mu$-dependent matrix $({A_\mu})_{i,j}=a_\mu(\theta_j,\theta_i)$, where $\left(\theta_i\right)_{1\leq i\leq n}$ denote the basis
functions of the
considered finite element space on $\Gamma$. Two different meshes, on which the matrices are assembled, are considered, see Table
\ref{tabmesh} and Figure \ref{mesh1}. The in-house code ACTIPOLE developed by EADS-IW and Airbus \cite{actipole1, actipole2} is used.
This test case is a challenging benchmark for two reasons.
First, the Green kernel $G_\mu(x,y):=\frac{\exp\left(i\mu\left|x-y\right|\right)}{4\pi\left|x-y\right|}$ oscillates
at a frequency proportional to the
parameter of interest $\mu$, and, secondly, the obtained matrices are dense and complex-valued.
Mesh~2 leads to a very large matrix and cannot be stored in an average desktop computer RAM.
The tests on Mesh~1 have been computed on a simple laptop with 4 Go of RAM, whereas the tests on Mesh~2 have been
computed on CCRT's Curie supercomputer \cite{curie}.

\renewcommand{\arraystretch}{1.5}
\begin{table}
\begin{center}
\begin{tabular}{|c |c |c |}
 \cline{2-3}
\multicolumn{1}{c|}{} & Mesh~1 & Mesh~2 \\
 \hline
number of triangles & $7,886$ & $40,576$ \\
 \hline
number of vertices  & $3,945$ & $20,290$ \\
 \hline
smallest edge (mm) & $6.53$ & $6.53$ \\
 \hline
mean edge (mm) & $437.52$ & $192.92$ \\
 \hline
largest edge (mm) & $718.99$ & $389.29$  \\
 \hline
number of complex nonzero coefficients per matrix & $1.56\times 10^7$ & $4.12\times 10^8$ \\
 \hline
memory usage to store one matrix in binary format (Go) & $0.23$ & $6.5$ \\
 \hline
\end{tabular}
\end{center}
\caption{\label{tabmesh} Characteristics of the two considered meshes.}
\end{table}

\begin{figure}[h!]
 \centering
\includegraphics[width=0.47\textwidth]{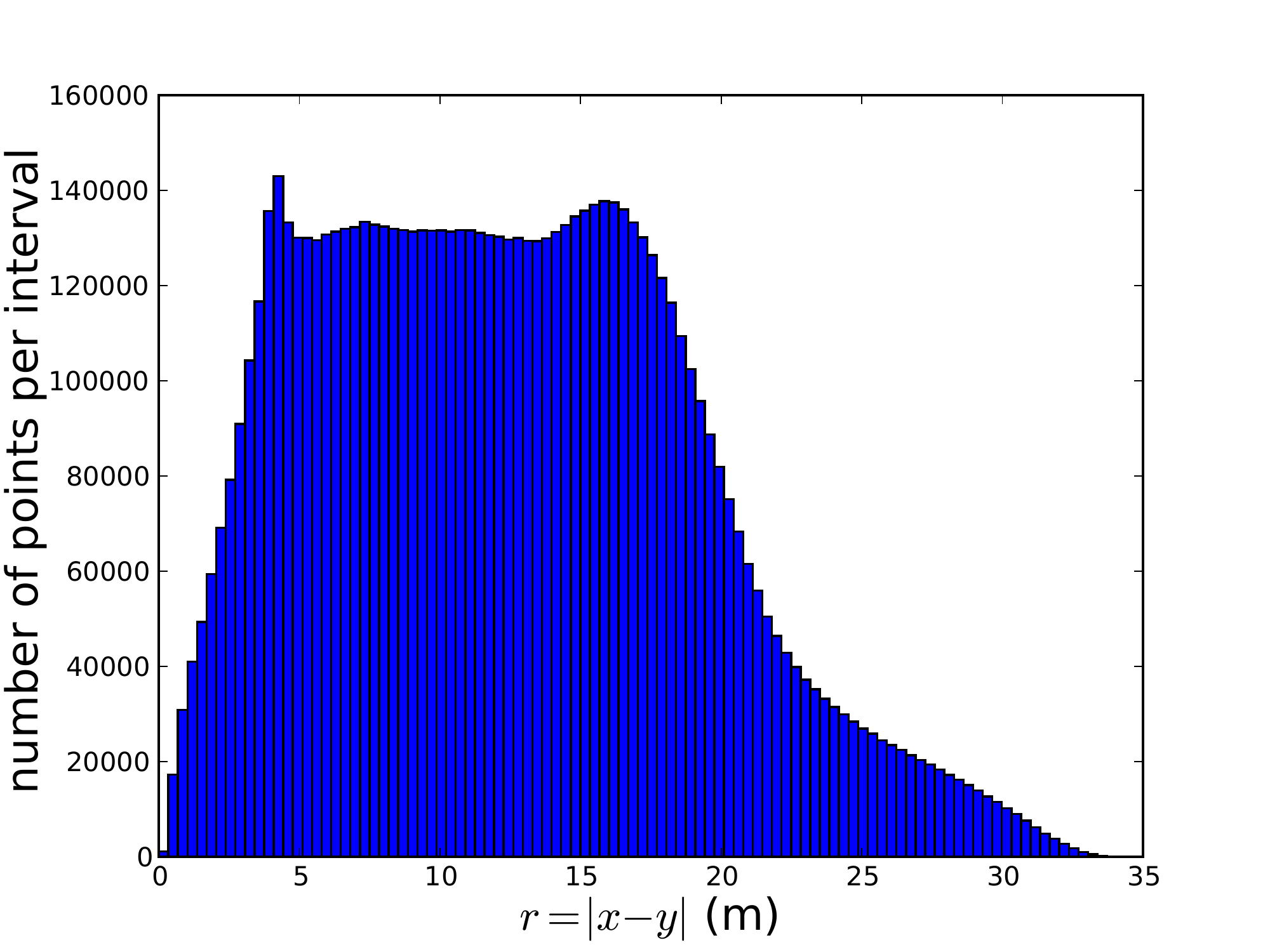}
\includegraphics[width=0.47\textwidth]{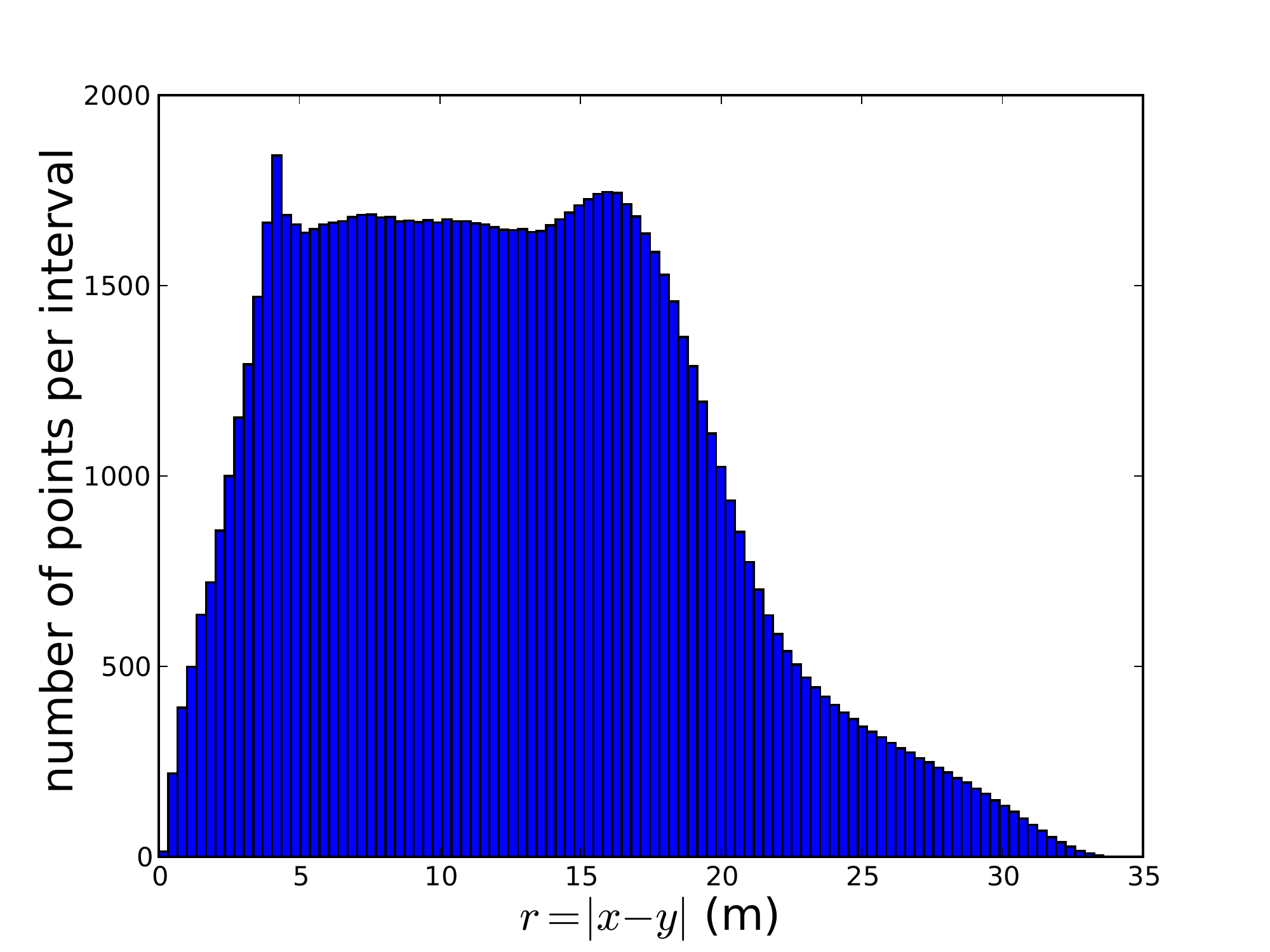}
 \caption{Histograms: discrete values of $r=|x-y|$ over the Gauss points from Mesh~1 (left), and chosen set of size $10^5$ (right).}
\label{fig:hist}
\end{figure}

To derive the approximation~\eqref{eq:non-intrusive} for $A_\mu$,
we carry out ${\rm EIM}^g$ to approximate
\begin{equation}
g(\mu,r):=\exp\left(i\mu r\right),~r=\left|x-y\right|,~x,y\in\Gamma.
\end{equation}
We choose $\mu\in\mathcal{P}_{\rm trial}:=\{0.005, 0.01, ..., 2.5\}$, a set of $1000$ values for the wave number,
so that the highest wave number for the source corresponds to a wavelength 5 times larger than the mean edge of Mesh~1.
A natural choice for the discrete set of values for $x$ and $y$ is the set
of Gauss points associated with the considered mesh, on which the quadrature formulae used to compute the integrals \eqref{eq:biliform}
are defined. The associated discrete set of values for $r=|x-y|$ is roughly proportional to the square of the
number of Gauss points, and equals $7.8\times 10^6$ for Mesh~1. To reduce the computational cost, a subsample of $10^5$
values for $r$, that has a very close density to the one obtained from the set of Gauss points, is chosen, see Figure \ref{fig:hist}.

Once ${\rm EIM}^g$ has been carried out, we can write
\begin{equation*}
A_\mu\approx \left(1+\mu^2\right)\sum_{m=1}^{d^g}\lambda^g_m(\mu)M_m,
\end{equation*}
where the matrices $M_m$ have been defined in Section~\ref{sec:descriptionproc}, so that the approximation~\eqref{eq:non-intrusive}
can be written using
\begin{equation}
\label{eq:XreuseEIMa}
z_p(\mu):=\left\{
\begin{alignedat}{2}
&\lambda^g_m(\mu),&\qquad &1\leq m\leq d^g, \quad p=m,\\
&\mu^2\lambda^g_m(\mu),&\qquad &1\leq m\leq d^g, \quad p=m+d^g.
\end{alignedat}\right.
\end{equation}
Note that we exploited the links in the functional dependence on $\mu$ for the two terms on the right-hand side of~\eqref{eq:biliform}
to carry out only one ${\rm EIM}^g$ procedure.

${\rm EIM}^g$ and ${\rm EIM}^z$ are carried out with respectively $d^g=30$ and ${d}^z=32$ interpolation points
(notice that $d_{\rm max}=60$).
To check the accuracy of the
approximation, we compute the relative Frobenius norm error on the matrix $A_\mu$ and the
relative Euclidian norm error on the acoustic pressure computed using
the approximate matrix, on a network of $400$ points located behind the aircraft. Figure~\ref{fig:diffmat} presents the results on Mesh~1.
In this figure, the relative differences are computed on $100$ values of $\mu$, namely one tenth of the considered parameter
values, explaining why only 7 minima are achieved on the left plot. On the right plot concerning the acoustic pressure behind the
aircraft, a large number of values are at the level of machine precision.
Note that the right-hand side of~\eqref{eq:hardsound} also depends on the parameter $\mu$. To compute the right plot of
Figure~\ref{fig:diffmat}, we computed the exact values of this right-hand.

\setlength\figureheight{0.43\textwidth}
\setlength\figurewidth{0.45\textwidth}
\begin{figure}[h!]
   \begin{minipage}[c]{.48\linewidth}
%
%
%
%
\begin{tikzpicture}

\begin{axis}[
xlabel={frequency (Hz)},
ylabel={relative error on the matrix},
xmin=0, xmax=140,
ymin=-16, ymax=-10,
axis on top,
width=\figurewidth,
height=\figureheight
]
\addplot [blue]
coordinates {
(0.270563403256222,-16) (1.62202624834986,-11.4877321463407) (2.97348909344351,-11.5105855029933) (4.32495193853715,-11.6336846194151) (5.67641478363079,-15.2235034049388) (7.02787762872443,-11.6829719360901) (8.37934047381807,-11.8855793818654) (9.73080331891171,-11.2704882553088) (11.0822661640054,-11.2011461925881) (12.433729009099,-11.2685891036591) (13.7851918541926,-11.4965620559469) (15.1366546992863,-11.4242838230464) (16.4881175443799,-11.4466367558778) (17.8395803894736,-11.3981514667148) (19.1910432345672,-11.8036391212847) (20.5425060796608,-11.9162317170956) (21.8939689247545,-15.739310453154) (23.2454317698481,-11.3809101926261) (24.5968946149418,-11.4156675915143) (25.9483574600354,-11.5513579309171) (27.2998203051291,-11.7580632606428) (28.6512831502227,-11.7065418746292) (30.0027459953163,-12.0287853731222) (31.35420884041,-15.6049430490138) (32.7056716855036,-11.7579549813393) (34.0571345305973,-11.6898289058567) (35.4085973756909,-11.841077570147) (36.7600602207846,-11.8304428801699) (38.1115230658782,-11.6601951013471) (39.4629859109718,-11.8038126012688) (40.8144487560655,-12.0706485845224) (42.1659116011591,-11.8782949518985) (43.5173744462528,-11.6201178062468) (44.8688372913464,-11.6318470606197) (46.22030013644,-11.744841208704) (47.5717629815337,-11.8713132435583) (48.9232258266273,-11.5473277824294) (50.274688671721,-11.4245781680749) (51.6261515168146,-11.7626630499233) (52.9776143619083,-11.7426257253766) (54.3290772070019,-11.3516091345522) (55.6805400520955,-11.2145653953587) (57.0320028971892,-11.3152727314305) (58.3834657422828,-11.7803947964709) (59.7349285873765,-11.3082950961466) (61.0863914324701,-11.0389903147196) (62.4378542775637,-10.8839036523815) (63.7893171226574,-10.8980041130612) (65.140779967751,-10.978512647994) (66.4922428128447,-11.251444215142) (67.8437056579383,-16) (69.195168503032,-11.3383916474497) (70.5466313481256,-11.1844105878902) (71.8980941932192,-11.2482745761718) (73.2495570383129,-11.5530378754885) (74.6010198834065,-12.2281710212722) (75.9524827285002,-11.46711242446) (77.3039455735938,-11.4060953460458) (78.6554084186874,-11.4871683990509) (80.0068712637811,-11.8502315388394) (81.3583341088747,-11.9000020221132) (82.7097969539684,-11.7623411543759) (84.061259799062,-11.7427506385594) (85.4127226441557,-12.057906492618) (86.7641854892493,-11.8371516372394) (88.1156483343429,-11.7298039205241) (89.4671111794366,-12.0699337038243) (90.8185740245302,-11.909931380482) (92.1700368696239,-11.6621607663204) (93.5214997147175,-11.4611125253832) (94.8729625598111,-11.8015022754312) (96.2244254049048,-12.2991792349715) (97.5758882499984,-11.5302332646964) (98.9273510950921,-11.2977930676096) (100.278813940186,-11.1657016145514) (101.630276785279,-11.1846825884283) (102.981739630373,-11.1890712812804) (104.333202475467,-11.25573188182) (105.68466532056,-11.1777037440344) (107.036128165654,-11.1778443856379) (108.387591010748,-11.0959494851753) (109.739053855841,-11.0600237772365) (111.090516700935,-11.1577482890883) (112.441979546028,-11.2861352456353) (113.793442391122,-11.4965281905113) (115.144905236216,-11.9394128905618) (116.496368081309,-11.6668164965238) (117.847830926403,-11.1942078032543) (119.199293771497,-10.9353551264197) (120.55075661659,-16) (121.902219461684,-11.2996250374509) (123.253682306778,-11.1219175957054) (124.605145151871,-11.1217396755515) (125.956607996965,-11.205605955569) (127.308070842059,-11.1121434020077) (128.659533687152,-11.1842804944179) (130.010996532246,-11.355499872152) (131.362459377339,-12.114864112936) (132.713922222433,-15.751933880241) (134.065385067527,-10.9980311354376) 
};
\path [draw=black, fill opacity=0] (axis cs:13,1)--(axis cs:13,1);

\path [draw=black, fill opacity=0] (axis cs:1,13)--(axis cs:1,13);

\path [draw=black, fill opacity=0] (axis cs:13,0)--(axis cs:13,0);

\path [draw=black, fill opacity=0] (axis cs:0,13)--(axis cs:0,13);

\end{axis}

\end{tikzpicture}
   \end{minipage} \hfill
   \begin{minipage}[c]{.48\linewidth}
%
%
%
%
\begin{tikzpicture}

\begin{axis}[
xlabel={frequency (Hz)},
ylabel={relative error on the solution},
xmin=0, xmax=140,
ymin=-16, ymax=-7,
axis on top,
width=\figurewidth,
height=\figureheight
]
\addplot [blue]
coordinates {
(0.270563403256222,-16) (1.62202624834986,-7.85130339869544) (2.97348909344351,-16) (4.32495193853715,-7.51932238058046) (5.67641478363079,-16) (7.02787762872443,-16) (8.37934047381807,-7.82102904745654) (9.73080331891171,-8.73983587358816) (11.0822661640054,-7.95057413156987) (12.433729009099,-16) (13.7851918541926,-16) (15.1366546992863,-9.09915643212437) (16.4881175443799,-16) (17.8395803894736,-16) (19.1910432345672,-11.6085408680016) (20.5425060796608,-16) (21.8939689247545,-16) (23.2454317698481,-16) (24.5968946149418,-16) (25.9483574600354,-16) (27.2998203051291,-16) (28.6512831502227,-16) (30.0027459953163,-16) (31.35420884041,-16) (32.7056716855036,-16) (34.0571345305973,-16) (35.4085973756909,-16) (36.7600602207846,-16) (38.1115230658782,-16) (39.4629859109718,-16) (40.8144487560655,-16) (42.1659116011591,-16) (43.5173744462528,-16) (44.8688372913464,-16) (46.22030013644,-7.45444631685031) (47.5717629815337,-16) (48.9232258266273,-16) (50.274688671721,-16) (51.6261515168146,-16) (52.9776143619083,-16) (54.3290772070019,-16) (55.6805400520955,-16) (57.0320028971892,-16) (58.3834657422828,-16) (59.7349285873765,-16) (61.0863914324701,-8.4351750784764) (62.4378542775637,-9.43223192436642) (63.7893171226574,-16) (65.140779967751,-10.4254693756186) (66.4922428128447,-16) (67.8437056579383,-16) (69.195168503032,-16) (70.5466313481256,-16) (71.8980941932192,-16) (73.2495570383129,-16) (74.6010198834065,-16) (75.9524827285002,-16) (77.3039455735938,-16) (78.6554084186874,-16) (80.0068712637811,-16) (81.3583341088747,-16) (82.7097969539684,-16) (84.061259799062,-16) (85.4127226441557,-16) (86.7641854892493,-16) (88.1156483343429,-11.4373666385081) (89.4671111794366,-16) (90.8185740245302,-16) (92.1700368696239,-16) (93.5214997147175,-16) (94.8729625598111,-11.4497263605173) (96.2244254049048,-16) (97.5758882499984,-16) (98.9273510950921,-16) (100.278813940186,-16) (101.630276785279,-16) (102.981739630373,-16) (104.333202475467,-16) (105.68466532056,-16) (107.036128165654,-16) (108.387591010748,-16) (109.739053855841,-16) (111.090516700935,-16) (112.441979546028,-16) (113.793442391122,-16) (115.144905236216,-16) (116.496368081309,-16) (117.847830926403,-16) (119.199293771497,-9.48977493785621) (120.55075661659,-16) (121.902219461684,-16) (123.253682306778,-7.33983203134236) (124.605145151871,-16) (125.956607996965,-16) (127.308070842059,-16) (128.659533687152,-16) (130.010996532246,-16) (131.362459377339,-16) (132.713922222433,-16) (134.065385067527,-16) 
};
\path [draw=black, fill opacity=0] (axis cs:13,1)--(axis cs:13,1);

\path [draw=black, fill opacity=0] (axis cs:1,13)--(axis cs:1,13);

\path [draw=black, fill opacity=0] (axis cs:13,0)--(axis cs:13,0);

\path [draw=black, fill opacity=0] (axis cs:0,13)--(axis cs:0,13);

\end{axis}

\end{tikzpicture}
   \end{minipage}
 \caption{Mesh~1, log${}_{10}$ of the relative error on the Frobenius norm of the matrix $A_\mu$ (left), and on the acoustic pressure computed using
the approximate matrix on a network of $400$ points located behind the plane in Euclidian norm (right), with $d^g=30$ and ${d}^z=32$.}
\label{fig:diffmat}
\end{figure}
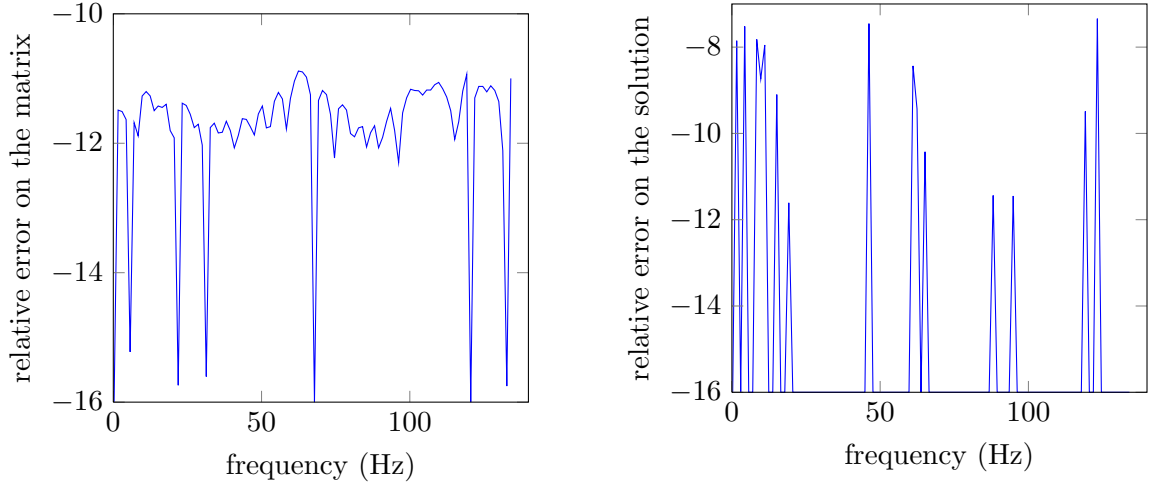

\begin{figure}[h!]
 \centering
\includegraphics[width=0.47\textwidth]{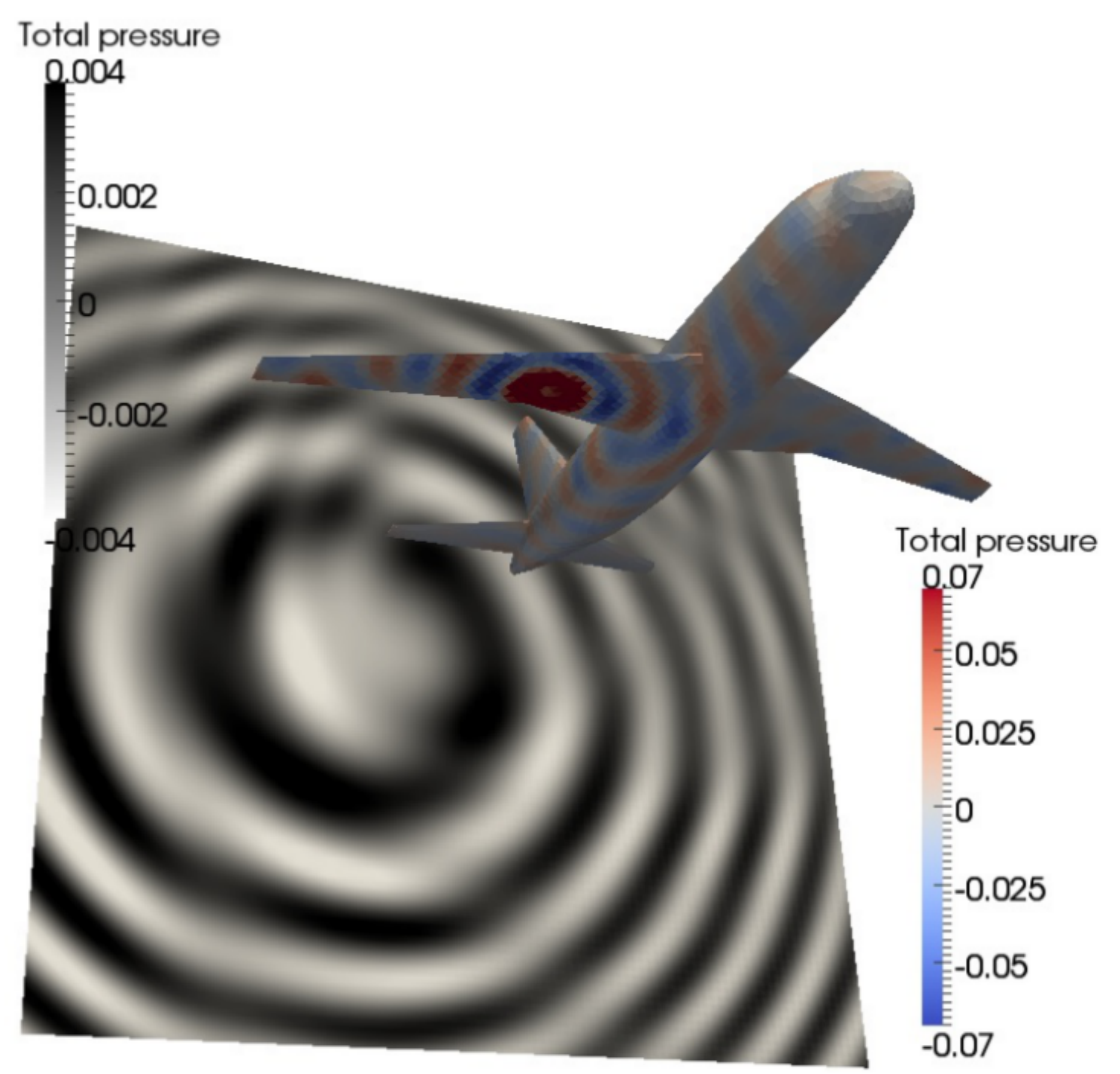}
\includegraphics[width=0.47\textwidth]{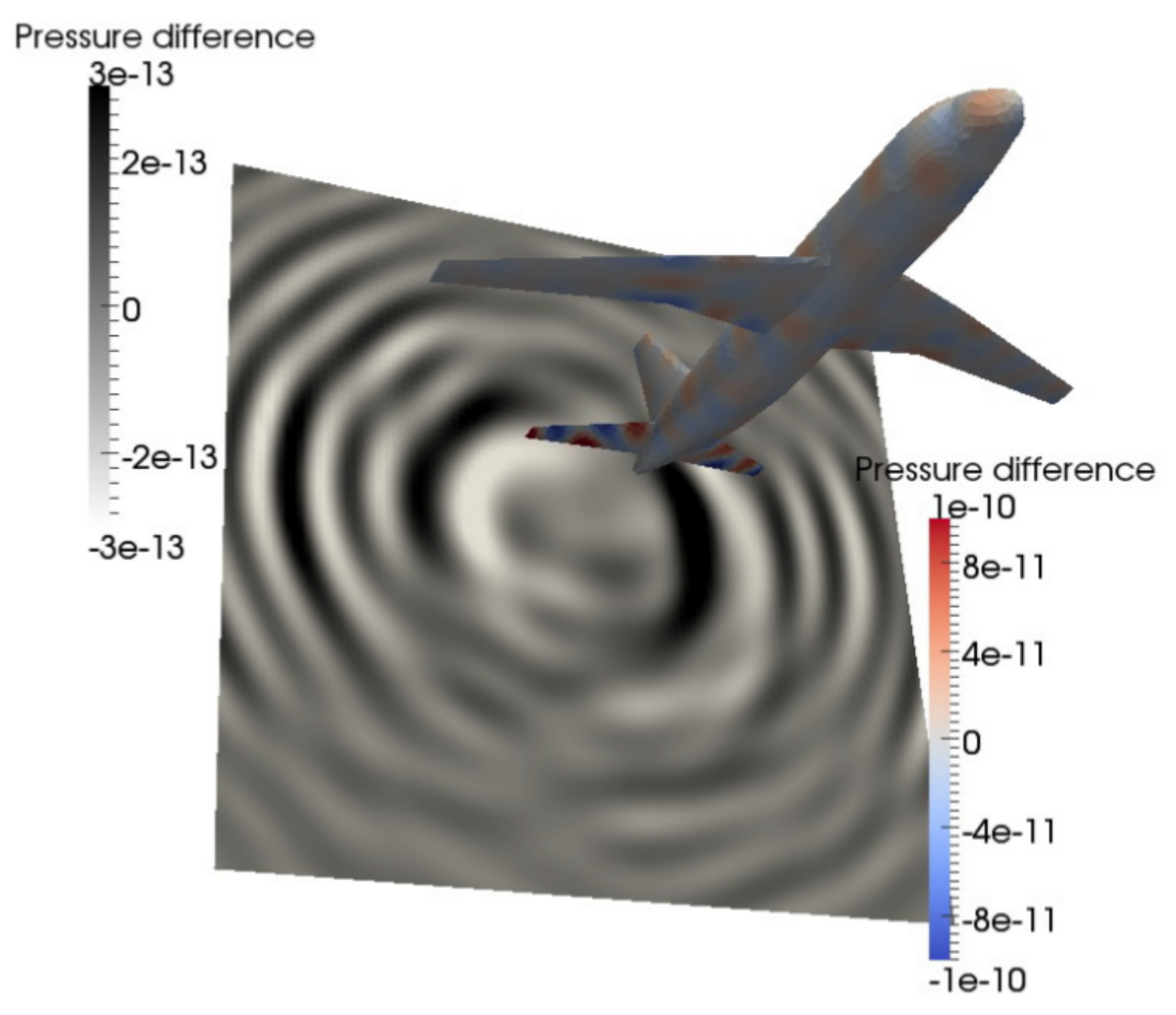}
 \caption{Mesh~1: total acoustic field on the plane and on the network of points (left), and difference between the
exact and approximate solution (right), for $\mu=2.47$.}
\label{fig:res1}
\end{figure}

Figure \ref{fig:res1} shows the solution to the problem on Mesh~1 and the relative difference of the solution using
the exact matrix and its approximation for $\mu=2.47$.

The simulation is repeated on Mesh~2, with $d^g=50$ and ${d}^z=50$. A twice as large frequency interval is considered
since Mesh~2 has a better spatial resolution than Mesh~1.
Figure \ref{fig:diffmat22} shows the relative Frobenius norm error on the matrix $A_\mu$, confirming the accuracy of the
approximation.

\setlength\figureheight{0.43\textwidth}
\setlength\figurewidth{0.45\textwidth}
\begin{figure}[h!]
\begin{centering}
%
%
%
%
\begin{tikzpicture}

\begin{axis}[
xlabel={frequency (Hz)},
ylabel={relative error on the matrix},
xmin=0, xmax=300,
ymin=-16, ymax=-11,
axis on top,
width=\figurewidth,
height=\figureheight
]
\addplot [blue]
coordinates {
(0.270563403256222,-16) (2.16450722604978,-12.2016507199209) (4.32901445209955,-12.0384872595277) (6.49352167814933,-11.6779683693134) (8.65802890419911,-11.4068451775321) (10.8225361302489,-11.1542822743988) (12.7164799530424,-11.3308583014645) (14.8809871790922,-11.5327774440492) (17.045494405142,-11.3986573670555) (19.2100016311918,-11.8374991496571) (21.3745088572415,-11.5970558895887) (23.2684526800351,-12.2969723344436) (25.4329599060849,-11.8741757584646) (27.5974671321347,-11.7914980852317) (29.7619743581844,-11.6695850935889) (31.9264815842342,-12.0207176008846) (34.090988810284,-11.8450170582418) (35.9849326330775,-12.2571466562295) (38.1494398591273,-12.3024501553717) (40.3139470851771,-11.5348826130437) (42.4784543112269,-11.9208876906617) (44.6429615372766,-11.7589441469243) (46.5369053600702,-11.6597860022347) (48.70141258612,-11.8082771564591) (50.8659198121698,-12.1080086811171) (53.0304270382195,-11.5664353363332) (54.3832440545006,-11.6620907200628) (55.1949342642693,-11.5699308792758) (57.0888780870629,-11.3787021816472) (59.2533853131126,-11.14975451322) (61.4178925391624,-11.1443368889534) (63.5823997652122,-11.1522179900464) (65.746906991262,-11.6317817185141) (67.9114142173117,-11.9868619180719) (69.8053580401053,-11.7422157891197) (71.9698652661551,-11.9970841170526) (74.1343724922048,-11.7823222128364) (76.2988797182546,-11.7960617535755) (78.4633869443044,-11.9138892869077) (80.3573307670979,-11.8812550083244) (82.5218379931477,-11.5090366925871) (84.6863452191975,-11.7081533438041) (86.8508524452473,-12.1746154131069) (89.0153596712971,-11.6273945764094) (90.3681766875782,-11.3846561830435) (90.9093034940906,-11.5144144651674) (93.0738107201404,-11.8365958138156) (95.2383179461902,-12.4214779750442) (97.4028251722399,-11.6287053934961) (99.5673323982897,-11.7021960072617) (101.731839624339,-11.3484954922608) (103.625783447133,-11.5634406116234) (105.790290673183,-11.9501932338398) (107.954797899233,-16) (110.119305125282,-11.707968260221) (112.283812351332,-12.1689679271971) (114.177756174126,-11.9233786057279) (116.342263400175,-11.4781945530865) (118.506770626225,-11.2542959957786) (120.671277852275,-11.5004820368989) (122.835785078325,-11.7127405712696) (124.729728901118,-11.6638907551572) (126.894236127168,-11.531939222874) (129.058743353218,-11.4860122903242) (131.223250579268,-11.8162466908706) (133.387757805317,-11.4115878100151) (135.552265031367,-11.4181586968835) (137.446208854161,-11.6516823500113) (139.610716080211,-11.4280500278254) (141.77522330626,-11.3277340140877) (143.93973053231,-11.0033232505127) (146.10423775836,-11.2591371535295) (147.998181581153,-11.3477636619352) (150.162688807203,-11.7038175111053) (152.327196033253,-11.6197286654693) (154.491703259303,-11.3240390197729) (156.656210485353,-11.2753660820934) (158.550154308146,-11.6904077683186) (160.714661534196,-16) (162.879168760246,-11.6646167725992) (165.043675986295,-11.4289627126626) (167.208183212345,-11.3216774399049) (169.372690438395,-11.4281563603795) (171.266634261189,-11.769486826976) (173.431141487238,-11.4966240797454) (175.595648713288,-11.383233152547) (177.760155939338,-11.4264485314354) (179.924663165388,-11.5544212642861) (180.4657899719,-11.6149601610225) (181.818606988181,-16) (183.983114214231,-11.6991957820553) (186.147621440281,-11.7173780322707) (188.312128666331,-11.7548097078778) (190.47663589238,-11.4968631060135) (192.370579715174,-11.5130625683851) (194.535086941224,-11.4149883203115) (196.699594167273,-11.4696844580669) (198.864101393323,-11.4045491480528) (201.028608619373,-11.4560462072073) (203.193115845423,-11.5810095161805) (205.087059668216,-11.6873184207321) (207.251566894266,-11.4473139069316) (209.416074120316,-11.2955112245507) (211.580581346366,-11.2169734735465) (213.745088572415,-11.0273604409499) (215.639032395209,-11.4106999987377) (219.968046847309,-11.4180115549237) (224.297061299408,-11.2125581027941) (228.355512348251,-11.3094540658417) (232.684526800351,-11.3390188409236) (237.013541252451,-11.3674989280828) (241.071992301294,-11.2696181837622) (245.401006753393,-11.265289210371) (249.459457802237,-12.0553008153918) (253.788472254336,-12.1570195937219) (258.117486706436,-11.3807855016996) (262.175937755279,-11.2742149616057) (266.504952207379,-11.1084674209745) (270.563403256222,-16) 
};
\path [draw=black, fill opacity=0] (axis cs:13,1)--(axis cs:13,1);

\path [draw=black, fill opacity=0] (axis cs:1,13)--(axis cs:1,13);

\path [draw=black, fill opacity=0] (axis cs:13,0)--(axis cs:13,0);

\path [draw=black, fill opacity=0] (axis cs:0,13)--(axis cs:0,13);

\end{axis}

\end{tikzpicture}
\end{centering}
 \caption{Mesh~2, log${}_{10}$ of the relative error on the Frobenius norm of the matrix $A_\mu$, with $d^g=50$ and ${d}^z=50$.}
\label{fig:diffmat22}
\end{figure}
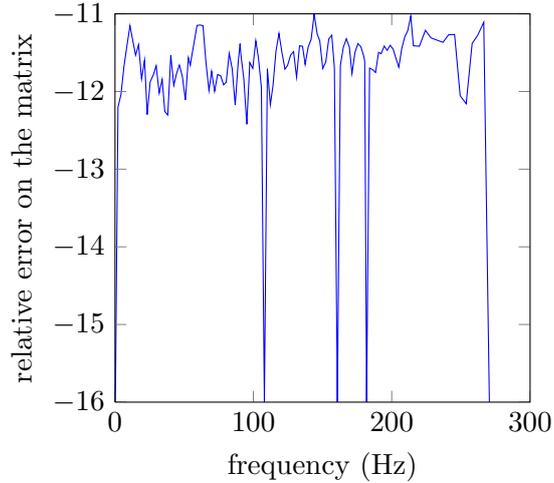

\subsection{Sound-hard scattering in a non-uniform flow}
\label{sec:application2}

\begin{figure}[h!]
 \centering
\includegraphics[width=0.47\textwidth]{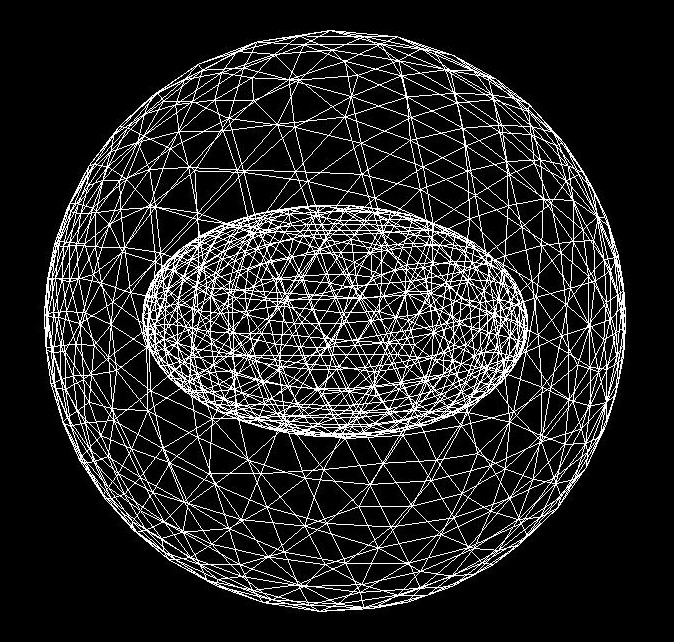}
 \caption{Representation of the mesh.}
\label{fig:mesh1}
\end{figure}

Consider an ellipsoid with major axis directed along the $z$-axis. This object is included inside a larger ball,
see Figure~\ref{fig:mesh1}.
The external border of the ball after discretization is denoted by $\Gamma_\infty$.
The complement of the ellipsoid in the ball is denoted by $\Omega^-$.
A potential flow is precomputed around the ellipsoid and inside the ball,
such that the flow is uniform outside the ball, of Mach number $0.3$ and directed along the $z$-axis.
The flow is fixed, and does not depend on the parameter $\mu$.
An acoustic monopole source lies upstream of the ball, on the $z$-axis as well.
The parameter is again the wave number of the monopole source.

The considered formulation is a coupled Finite Element Method (FEM) - BEM formulation described in \cite{jcp}.
It consists in (i) applying a change of variable to transform the convected Helmholtz equation into the classical Helmholtz equation
outside the ball, in order to apply a standard BEM, and (ii) stabilizing the formulation
to avoid resonant frequencies associated with the eigenvalues of the Laplacian inside the ball of border $\Gamma_\infty$.
The formulation depends on the wave number of the source in a complex way,
but we will see in our numerical tests that our nonintrusive procedure provides an accurate approximation of the resulting matrix as
a linear combination of a few snapshots of the complete matrix at some wave numbers of the source.

Consider the product space $\mathbb{H}:=H^{1}\left({\Omega^-}\right)\times H^{-\frac{1}{2}}\left({\Gamma_\infty}\right)\times H^{1}({\Gamma_\infty})$
with inner product
$\left(\left(\Phi,\lambda,p\right),\left(\Phi^t,\lambda^t,p^t\right)\right)_{\mathbb{H}}:=\left(\Phi, \Phi^t\right)_{H^{1}\left({\Omega^-}\right)}+\left(\lambda,\lambda^t\right)_{H^{-\frac{1}{2}}\left({\Gamma_\infty}\right)}
+\left(p,p^t\right)_{H^{1}\left({\Gamma_\infty}\right)}$.
The weak formulation is:
Find $\left(\Phi,\lambda,p\right)\in \mathbb{H}$ such that $\forall\left(\Phi^t,\lambda^t, p^t\right)\in \mathbb{H}$,
\begin{subequations}
\label{eq:weakcouptrans2}
\begin{align}
\mathcal{V}_\mu(\Phi,\Phi^t)+\left(N_\mu({\gamma_0^-}\Phi),{\gamma_0^-}\Phi^t\right)_{\Gamma_{\infty}}+\left(\left({\tilde{D}_\mu}-\frac{1}{2}I\right)(\lambda),{\gamma_0^-}\Phi^t\right)_{\Gamma_{\infty}}
&= \left({\gamma_1}{f_{\rm inc}}_\mu,{\gamma_0^-}\Phi^t\right)_{\Gamma_{\infty}},\label{eq:weakcouptrans21}\\
\left(\lambda^t,\left({D}_\mu-\frac{1}{2}I\right)({\gamma_0^-}\Phi)\right)_{\Gamma_{\infty}}-\left(\lambda^t,S_\mu(\lambda)\right)_{\Gamma_{\infty}}
-i\left(\lambda^t,p\right)_{\Gamma_{\infty}}
&= -\left(\lambda^t,{\gamma_0}{f_{\rm inc}}_\mu\right)_{\Gamma_{\infty}},\label{eq:weakcouptrans22}\\
\left(N_\mu({\gamma_0^-}\Phi),p^t\right)_{\Gamma_{\infty}}+\left(\left({\tilde{D}}_\mu+\frac{1}{2}I\right)(\lambda),p^t\right)_{\Gamma_{\infty}}
-\delta_{\Gamma_\infty}(p,p^t)&= \left({\gamma_1}{f_{\rm inc}}_\mu,p^t\right)_{\Gamma_\infty},\label{eq:weakcouptrans23}
\end{align}
\end{subequations}
where $\left(\cdot,\cdot\right)_{\Gamma_\infty}$ denotes the extension of the $L^2(\Gamma_\infty)$-inner product to the duality
pairing on $H^{-\frac{1}{2}}\left(\Gamma_\infty\right)\times H^{\frac{1}{2}}\left(\Gamma_\infty\right)$, and where
\begin{equation}
\label{eq:defc}
\delta_{\Gamma_\infty}(p,q):=\left(\vec{\nabla}_{\Gamma_\infty}p, \vec{\nabla}_{\Gamma_\infty}q\right)_{\Gamma_\infty}+
\left(p,q\right)_{\Gamma_\infty},
\end{equation}
with $\vec{\nabla}_{\Gamma_\infty}$ the surfacic gradient on $\Gamma_\infty$, and
\begin{equation}
\mathcal{V}_\mu(\Phi,\Phi^t):=\int_{{\Omega^-}}\Xi\vec{\nabla}\overline{\Phi}\cdot\vec{\nabla}{{\Phi}^t}-\mu^2\int_{{\Omega^-}}
\beta \overline{\Phi}{{\Phi}^t}+i\mu \int_{{\Omega^-}}\vec{V}\cdot\left(\overline{\Phi}\vec{\nabla}{{\Phi}^t}-{{\Phi}^t}
\vec{\nabla}\overline{\Phi}\right),
\end{equation}
where 
$\beta := r\left(\left(\varsigma+\gamma_\infty^2 P\right)^2-\gamma_\infty^4 M_\infty^2\right)$,
$\vec{V} := r\left(\left(\varsigma+\gamma_\infty^2 P\right)\mathcal{N}\vec{M} -\gamma_\infty^3\vec{M}_\infty\right)$,
$\Xi := r\mathcal{N}\mathcal{O}\mathcal{N}$ with
$r := \frac{\rho}{\rho_\infty}$,
$\varsigma :=\frac{c_\infty}{c}$,
$\gamma_\infty := \frac{1}{\sqrt{1-M_\infty^2}}$,
$P := \vec{M}\cdot\vec{M}_\infty$,
$\mathcal{N} := I+C_\infty\vec{M}_\infty\vec{M}^T_\infty$,
$\mathcal{O} := I-\vec{M}\vec{M}^T$,
and $C_\infty := \frac{\gamma_\infty-1}{M_{\infty}^2}$.
In the above notation, the subscript $\infty$ is used for quantities outside the ball,
$\rho$ is the density of the flow, $c$ is the speed of sound when the flow is at rest and
$\vec{M}=\frac{\vec{v}}{c}$, where $\vec{v}$ is the velocity of the flow.
The operators $\gamma_0$ and $\gamma_1$ are Dirichlet and Neumann traces
on the coupling surface $\Gamma_\infty$.
The operators $N_\mu$, $D_\mu$, $\tilde{D}_\mu$, and $S_\mu$ are boundary integral operators, expressed in terms of
the Green kernel $G_\mu(x,y)=\frac{\exp(i\mu|x-y|)}{4\pi|x-y|}$ associated with the Helmholtz equation at wave number $\mu$.

The next step is to identify the dependencies in $\mu$ in the formulation~\eqref{eq:weakcouptrans2}.
It turns out that the functions of $\mu$ involved in the integrals of the formulation~\eqref{eq:weakcouptrans2} are
$\mu$, $\mu^2$, $\exp(i\mu r)$, $\mu \exp(i\mu r)$, $\mu^2 \exp(i\mu r)$, and $\mu\left(\frac{2i\pi\mu}{c}\mu-1\right)\exp(i\mu r)$.
As in the previous test case, ${\rm EIM}^g$ is carried out to approximate the function
$g(\mu,r)=\exp(i\mu r)$, $r=|x-y|$, $x,y\in\Gamma_\infty$.
We choose $\mu\in\mathcal{P}_{\rm trial}:=\{10, 10.03, ..., 40\}$, a set of $1000$ values for the wave number,
so that the highest wave number of the source corresponds to a wavelength 5 times larger than the mean edge of the mesh.
This time, instead of considering a subset of $|x-y|$ where $x$ and $y$ are the Gauss points associated with the mesh,
we take $r\in\{0,h,...,~Nh\}$, where $N=10000$ and $h=\frac{D}{N}$, $D$ being the diameter of the sphere $\Gamma_\infty$.
With this choice, we no longer need to know the position of the Gauss points, but simply the diameter of the
geometry of the test case.

Then, the functions $\left(\lambda^g_m(\mu)\right)_{1\leq m\leq d^g}, \mu\in\mathcal{P}_{\rm trial}$, are computed using
\eqref{eq:onlinea1pb}, and the functions $z_p(\mu)$, $1\leq p\leq d_{\rm max}:=4d^g+3$, are defined by
\begin{equation}
z_p(\mu):=\left\{
\begin{alignedat}{4}
&\lambda^g_m(\mu),&\qquad 1&\leq p\leq d^g,&\quad &m=p,\\
&\mu\lambda^g_m(\mu),&\qquad d^g+1&\leq p\leq 2d^g,&\quad &m=p-d^g,\\
&\mu^2\lambda^g_m(\mu),&\qquad 2d^g+1&\leq p\leq 3d^g,&\quad &m=p-2d^g,\\
&\mu\left(\frac{2i\pi}{c}\mu-1\right)\lambda^g_m(\mu),&\qquad 3d^g+1&\leq p\leq 4d^g,&\quad &m=p-3d^g,\\
&1,&\qquad  p&=4d^g+1,&&\\
&\mu,&\qquad  p&=4d^g+2,&&\\
&\mu^2,&\qquad  p&=4d^g+3.&&
\end{alignedat}\right.
\end{equation}
${\rm EIM}^g$ and ${\rm EIM}^z$ are carried out with respectively $17$ and $20$ interpolation points (notice that $d_{\rm max}=71$).

\begin{figure}[h!]
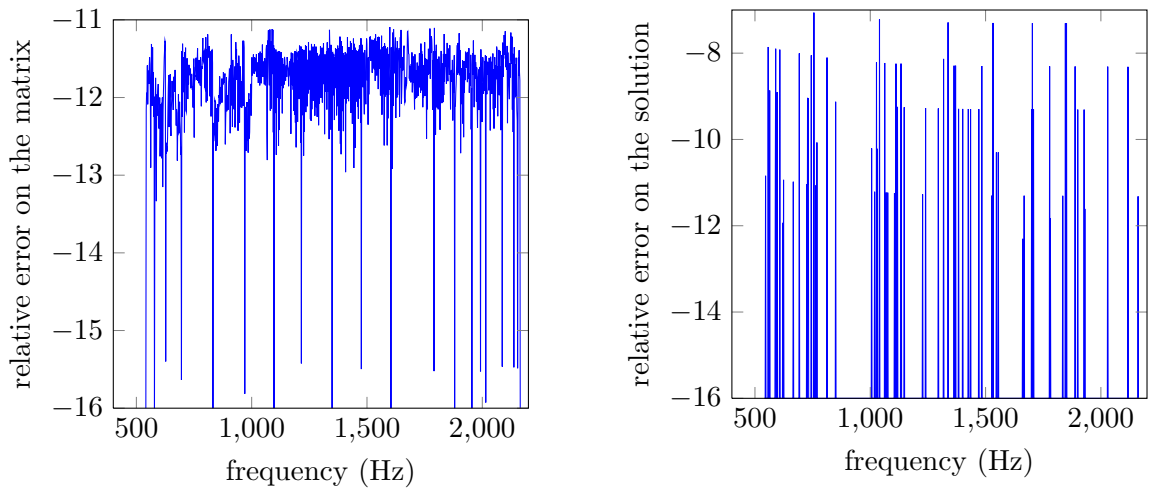

   \begin{minipage}[c]{.48\linewidth}
\include{erreurmatred2.tikz}
   \end{minipage} \hfill
   \begin{minipage}[c]{.48\linewidth}
\include{erreur_proche2.tikz}
   \end{minipage}
 \caption{Log${}_{10}$ of the relative error on the Frobenius norm of the matrix $A_\mu$ (left) and on the acoustic pressure computed using
the approximate matrix computed using \eqref{eq:non-intrusive} on a network of $400$ points located behind the object in Euclidian norm (right), with $d^g=17$ and ${d}^z=20$.}
\label{fig:diffmat2}
\end{figure}

Figure \ref{fig:diffmat2} shows
the relative Frobenius norm error on the matrix $A_\mu$ and the relative Euclidian norm error on the acoustic pressure computed using
the approximate matrix on a network of $400$ points located behind the scattering ellipsoid.
In this test case, an excellent accuracy is obtained with only $20$ precomputed matrices.

\setlength\figureheight{0.43\textwidth}
\setlength\figurewidth{0.45\textwidth}

\section{Conclusion and outlook}
\label{sec:otherapprox}
The method described herein provides an efficient nonintrusive approximation of parameter-dependent linear
systems, provided that the considered code can return the assembled matrix and that the corresponding weak formulation is known.
The method offers a crucial practical advantage over existing methods since it avoids significant implementation efforts.
In the present work, the choice has been made to approximate the whole matrix $A_\mu$ assembled by the code, but the procedure applies in the same way to
the approximation of any linear functional $l$ of the matrix $A_\mu$, whereby
\begin{equation}
\label{eq:non-intrusive_lin}
l(A_\mu)\approx \sum_{m=1}^{{d}^z} \beta^z_{m}(\mu) l(A_{\mu^z_m}),
\end{equation}
where the storage of $A_{\mu^z_m}$ for all $1\leq m\leq {{d}^z}$ is replaced by the storage of $l(A_{\mu^z_m})$ for all
$1\leq m\leq {{d}^z}$, which may be much lighter in terms of memory usage.
The efficient construction of the reduced matrix $\hat{A}_\mu$ in the RBM corresponds to $l(A_\mu)=U^t A_\mu U$, as explained in the
introduction.

Finally, we observe that in the case where the right-hand side $C$ of the problem \eqref{eq:matrix_pb} is also dependent on the
parameter $\mu$ (then written $C_\mu$), the same procedure can be applied to derive a separated representation of $C_\mu$.

\section*{Acknowledgement}
This work was partially supported by EADS Innovation Works.


\begin{thebibliography}{1}

\bibitem{curie}
http://www-hpc.cea.fr/en/complexe/tgcc-curie.htm.

\bibitem{RBconv}
P.~Binev, A.~Cohen, W.~Dahmen, R.~A. DeVore, G.~Petrova, and P.~Wojtaszczyk.
\newblock Convergence rates for greedy algorithms in reduced basis methods.
\newblock {\em SIAM J. Math. Analysis}, pages 1457--1472, 2011.

\bibitem{jcp}
F.~Casenave, A.~Ern, and G.~Sylvand.
\newblock A coupled boundary element/finite element method for the convected
  {H}elmholtz equation with non-uniform flow in a bounded domain.
\newblock {\em arXiv preprint arXiv:1303.6923}, 2013.

\bibitem{actipole1}
A.~Delnevo and I.~Terrasse.
\newblock Code {ACTI3S} harmonique, justification math\'ematique, {P}artie {I}.
\newblock Technical report, EADS, 2001.

\bibitem{actipole2}
A.~Delnevo and I.~Terrasse.
\newblock Code {ACTI3S}, justifications math\'ematiques, {P}artie {II} :
  presence d'un \'ecoulement uniforme.
\newblock Technical report, EADS, 2002.

\bibitem{maday}
Y.~Maday, N.C. Nguyen, A.T. Patera, and S.~Pau.
\newblock A general multipurpose interpolation procedure: the magic points.
\newblock {\em Communications On Pure And Applied Analysis}, 8(1):383--404,
  2008.

\bibitem{nedelecbook}
J.C. N\'ed\'elec.
\newblock {\em Acoustic and Electromagnetic Equations: Integral Representations
  for Harmonic Problems}.
\newblock Number vol.~144 in Applied Mathematical Sciences. Springer, 2001.

\bibitem{RB}
C.~Prud'homme, D.V. Rovas, K.~Veroy, L.~Machiels, Y.~Maday, A.T. Patera, and
  G.~Turinici.
\newblock Reliable real-time solution of parametrized partial differential
  equations: Reduced-basis output bound methods.
\newblock {\em CJ Fluids Engineering}, 124:70--80, 2002.

\bibitem{krig}
M.L. Stein.
\newblock {\em Interpolation of spatial data: some theory for kriging}.
\newblock Springer Series in Statistics Series. Springer London, Limited, 1999.

\end{thebibliography}
\end{document}